\documentclass[twocolumn,10pt, journal]{IEEEtran}

\usepackage{tikz}
\usepackage{amsthm}
\usepackage{bm, bbm}
\usepackage{blindtext}
\usepackage{graphicx}
\usepackage{amsmath}
\usepackage{mathrsfs}
\usepackage[capitalise]{cleveref}
\usepackage[sans]{dsfont}
\usepackage{amsbsy}
\usepackage{amssymb}
\usepackage{pgfplots}
\usetikzlibrary{arrows,positioning,patterns}

\newtheorem{assumption}{Assumption}
\newtheorem{problem}{Problem}

\newtheorem{remark}{Remark}
\newtheorem{theorem}{Theorem}
\newtheorem{definition}{Definition}

\newtheorem{example}{Example}

\newcommand{\Eqdef}{{:=}}

\newcommand{\range}{\! : \! }

\newcommand{\real}{\ensuremath{\mathbb{R}}}

\newcommand{\norm}[1]{\left\lVert#1\right\rVert}

\DeclareMathOperator{\EX}{\mathbb{E}}

\begin{document}

\title{Channels, Remote Estimation and Queueing Systems With A Utilization-Dependent Component: A Unifying Survey Of Recent Results}

\author{Varun Jog, Richard J. La,  
    Michael Lin, and Nuno C. Martins
\thanks{The authors are listed according to the lexicographical order of their last names. Varun Jog is with the Electrical and Computer Engineering Department at UW-Madison. Richard J. La, Michael Lin, and Nuno Miguel Lara Cintra Martins are with the Department of Electrical and Computer Engineering and the Institute for Systems Research, University of Maryland, College Park, MD, 20742 USA. E-mails: \texttt{vjog@wisc.edu, \{hyongla, mlin1025, nmartins\}@umd.edu}. Work writing this article was supported by AFOSR grant FA95501510367, NSF grant ECCS 1446785 and NIST grant 70NANB16H024}}

\maketitle

\begin{abstract}  

In this article, we survey the main models, techniques, concepts, and results centered on the design and performance evaluation of engineered systems that rely on a utilization-dependent component (UDC) whose operation may depend on its usage history or assigned workload.  Specifically, we report on research themes concentrating on the characterization of the capacity of channels and the design with performance guarantees of remote estimation and queueing systems. Causes for the dependency of a UDC on past utilization include the use of replenishable energy sources to power the transmission of information among the sub-components of a networked system, 
and the assistance of a human operator for servicing a queue. Our analysis unveils the similarity of the UDC models typically adopted in each of the research themes, and it reveals the differences in the objectives and technical approaches employed. We also identify new challenges and future research directions inspired by the cross-pollination among the central concepts, techniques and problem formulations of the research themes discussed.
\end{abstract}


\section{Introduction}
\label{sec:Introduction}
As new technologies and applications emerge, the algorithms that determine the functionality and regulate the operation of engineered systems have to contend with unexampled constraints and nonstandard problems. This evolution has been evident for communication~\cite{Ulukus2015Energy-harvesti}, cyber-physical~\cite{Kim2012Cyber-physical-,Baheti2011The-impact-of-c}, human-assisted~\cite{Staal2004Stress-cognitio}, networked estimation~\cite{Leong2018Optimal-control} and control~\cite{Hespanha2007A-survey-of-rec} systems, which are now designed for maximal performance subject to restrictions that are more intricate than the conventional limits on reliability and power usage. In this article, we provide a partial account of such advances by surveying models, concepts and results on the characterization of channel capacity, and the design and performance analysis of queuing and remote estimation systems, all of which have in common the unconventional attribute of relying on a component whose performance may be constrained by its usage history and possibly also be affected by the workload assigned to it.  We refer succinctly to this class of components as {\bf UDC}, which stands for \underline{utilization-dependent component}. We do not aim at a comprehensive survey; instead, we will cite a selection of published work relevant to each key concept, problem formulation or technique on an as-needed basis for illustration.

A primary goal of this article is to highlight the 
commonalities in the presented studies and describe
some of existing tools available to investigate 
challenging problems in these and other related fields. 
Our hope is that this survey article will serve as a 
good starting point for those who are interested
in conducting research in these areas and foster future 
research that builds on the cross-pollination among 
the methods and problem formulations originally developed 
and employed on each of the research themes broached.

\textbf{Paper structure:} After the Introduction, in Section~\ref{sec:UDCModelConstraints} we define a class of UDCs that is general enough to model the performance restrictions imposed by the reliance on the energy harvested from stochastic sources, human-assisted decision-making, or human labor. In Section~\ref{sec:CommonPOCMC}, we introduce widely-used models quantifying certain performance-limiting factors, such as mental workload, queueing workload and the state of charge of the battery of an energy harvesting module. Subsequently, in Sections~\ref{sec:Channels}~-~\ref{sec:Queues}, we employ these definitions as a unifying framework to discuss research on methods to design and analyze the performance of systems comprising a UDC in the context of communication, remote estimation, and queueing respectively. In Section~\ref{sec:UtilizationDependent}, we discuss recent work that addresses a class of problems that combines models, concepts and metrics from Sections~\ref{sec:RemoteEstimation} and~\ref{sec:Queues}, and, hence, serves as an example of the sort of research this survey intends to foster, and the UDC framework may facilitate. This article ends with the conclusions and future directions proposed in Section~\ref{sec:ConcFutDir}.

\textbf{Notation} Throughout the article, 
we use $\mathbb{N}$ and $\mathbb{N}_+$ to
denote the set of non-negative integers and
positive integers, respectively. When appropriate, 
we assume that all random 
variables and stochastic processes are defined
on a common probability space. In addition, 
for notational simplicity, we use 
$\mathcal{P}_{R_1 | R_2}(r_1 | r_2)$ to denote the
conditional probability ${\bf P}(R_1 = r_1
| R_2 = r_2)$, where $R_1, R_2$ are random
variables/vectors defined on the same probability
space. 

\section{A General Utilization-Dependent Component (UDC) Model}
\label{sec:UDCModelConstraints}

We start by presenting a model that is general enough to describe all types of UDC considered throughout this article. Without loss of generality, we limit our discussion to discrete-time processes and models. Namely, time takes values in the set of non-negative integers $\mathbb{N}$.

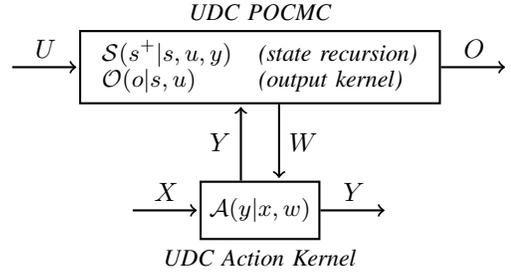
\begin{figure}
\centering

\begin{tikzpicture}[shorten >=1pt,on grid,auto]
    \tikzstyle{block}=[shape=rectangle,thick,draw,minimum size=0.3 in]

\node[block] (POCMC) {\small $\begin{array}{ll} \mathcal{S}(s^+|s,u,y) & \text{\it (state recursion)} \\ \mathcal{O}(o|s,u) & \text{\it (output kernel)} \end{array}$};
\node[anchor=south] at (POCMC.north) {\it \small UDC POCMC};

\node[below=0.4 in, block] at (POCMC.south) (ActionKernel) {\small $\mathcal{A}(y|x,w)$};
\node[anchor=north] at (ActionKernel.south) {\it \small UDC Action Kernel};


\draw[->,thick] (POCMC.east) -- ([xshift=0.35in]POCMC.east) node[midway,above] {$O$};

\draw[->,thick] ([xshift=-0.35in]POCMC.west) -- (POCMC.west) node[midway,above] {$U$};

\draw[->,thick] ([xshift=0.1in]POCMC.south) -- ([xshift=0.1in]ActionKernel.north) node[midway,right] {$W$}; 

\draw[->,thick] ([xshift=-0.1in]ActionKernel.north) -- ([xshift=-0.1in]POCMC.south) node[midway,left] {$Y$};

\draw[->,thick] (ActionKernel.east) -- ([xshift=0.35in]ActionKernel.east) node[midway,above] {$Y$};

\draw[->,thick] ([xshift=-0.35in]ActionKernel.west) -- (ActionKernel.west) node[midway,above] {$X$};

\end{tikzpicture}
\caption{Basic overall structure of the utilization-dependent component (UDC) model.}
\label{fig:UDCModel}
\end{figure}

\begin{definition}{\bf (UDC Model)}
The following describes the two main sub-components of the UDC model (see~Fig.~\ref{fig:UDCModel}).
\begin{itemize}
	\item The first sub-component is a partially-observed controlled Markov chain 
	(POCMC),\footnote{Here, we
	consider a POCMC in the model, rather than 
	a controlled Markov chain (CMC), to allow
	scenarios where the states are partially 
	observed. However, 
	this assumption is not critical for the purpose
	of this article.} whose state is represented as $S:=\{ S(k) : k \in \mathbb{N} \}$. It has two inputs denoted as $Y:=\{ Y(k) : k \in \mathbb{N} \}$ and $U:=\{ U(k) : k \in \mathbb{N} \}$, where the latter is an external control process. The outputs are indicated as $O:=\{ O(k) : k \in \mathbb{N} \}$ and $W:=\{ W(k) : k \in \mathbb{N} \}$. The former is characterized by an output kernel and is available to the policy that generates $U$, while the latter is a deterministic function of $S$ that we refer to as the performance process. 
	The processes $U$, $Y$, $S$ and $W$ take values in given alphabets $\mathbb{U}$, $\mathbb{Y}$, $\mathbb{S}$ and $\mathbb{W}$, respectively, which are subsets of real coordinate spaces. The POCMC is specified by maps $\mathcal{S}:\mathbb{S}^2 \times \mathbb{U} \times \mathbb{Y} \rightarrow [0,1]$, $\mathcal{O}:\mathbb{O} \times \mathbb{S} \times \mathbb{U}  \rightarrow [0,1]$ and $\mathcal{W}:\mathbb{U} \times \mathbb{S} \rightarrow \mathbb{W}$. The first two determine the state transition probability and the output kernel as follows:
	\begin{multline*}
	\mathcal{S}(s^+|s,u,y):=\mathcal{P}_{S(k+1)|S(k),U(k),Y(k)}(s^+|s,u,y), \\ \quad s,s^+ \in \mathbb{S}, \ u \in \mathbb{U}, \ y \in \mathbb{Y} 
	\end{multline*} 
	\begin{multline*}
	\mathcal{O}(o|s,u):=\mathcal{P}_{O(k)|S(k),U(k)}(o|s,u), \\ \quad o \in \mathbb{O}, \ s \in \mathbb{S}, \ u \in \mathbb{U}
	\end{multline*} The performance process is determined as ${ \mathcal{W}:\big( U(k),S(k) \big) \mapsto W(k)}$.
	\item The second sub-component is an action kernel that models the functionality whose performance is affected by the process $W$. Specifically, the output of the action kernel is $Y$ and the inputs are $W$ and an external source or command signal denoted as $X:=\{ X(k) : k \in \mathbb{N} \}$. The processes $X$ and $Y$ take values in given alphabets $\mathbb{X}$ and $\mathbb{Y}$ , respectively, which are subsets of real coordinate spaces. A map $\mathcal{A}:\mathbb{X}\times \mathbb{W} \rightarrow [0,1]$ specifies probabilistically $Y$ in terms of $X$ and $W$ as follows:
	\begin{multline*} \mathcal{A}(y|x,w):=\mathcal{P}_{Y(k)|X(k),W(k)}(y|x,w), \\ \quad y \in \mathbb{Y}, \ x \in \mathbb{X}, \ w \in \mathbb{W} 
	\end{multline*}
	Since $X$ and $U$ are both inputs of the overall UDC model, we allow causal policies that determine them in terms of~$O$. 
\end{itemize} 
The definition of the UDC model is not complete until we specify the probabilistic dependence among the implicit sources of randomness of the state recursion and the output kernels. These will be particularized throughout the text on an as-needed basis. Typically,  $S(k+1)$ and $O(k)$ are conditionally independent given $S(k)$, $U(k)$ and $Y(k)$; and $Y(k)$ and $O(k)$ are conditionally independent given $S(k)$, $U(k)$ and $X(k)$. 
\end{definition}


Although most existing work adopts variations of the 
models discussed in Section~\ref{sec:CommonPOCMC}, we 
opted to define a UDC model that is general enough to 
be used as a common framework to facilitate our
presentation in this article and also for future work.

\section{Commonly used POCMC models}
\label{sec:CommonPOCMC}

We proceed with defining a few common POCMC models, which we will invoke later appropriately altered to suit a specific application. In Sections~\ref{subsec:WorkloadUtilizationRatio} and~\ref{subsec:EH_Models}, we will specify $W$ as a function of $U$ and we will explicitly describe a state recursion in cases when it is clearer to do so. Because the output kernel is application dependent, we will defer its specification on an as-needed basis to Sections~\ref{sec:Channels}-\ref{sec:RemoteEstimation}. Typical cases include when $O$ equals $S$ or when additive measurement noise is present.

To be concise, when the POCMC is deterministic, we specify it via the functional recursion that governs the state update in terms of the inputs, and we express each output as a function of the current state and inputs. In the stochastic case, the probabilistic state recursion and output kernel can always be specified by the conditional probabilities associated with $\mathcal{S}$ and $\mathcal{O}$, respectively. The action kernel is specified in an analogous manner.

In order to appropriately indicate the dependence on certain parameters, or to discern which model is associated with a given internal process, such as $S$ and $W$, we often annotate them with a self-descriptive superscript.

\subsection{Utilization ratio and workload models}
\label{subsec:WorkloadUtilizationRatio} 

The concept of {\it mental workload}~\cite{Staal2004Stress-cognitio}, in human-assisted systems, refers broadly to the burden imposed on a human operator by the difficulty of tasks and the frequency with which they are assigned. Hence, considering that it is known to influence the performance of a human operator~\cite{Teigen94}, quantifying workload is important for the design of task assignment policies. In spite of having a rather simple structure, \cite[Chapter~11]{Wickens2000Engineering-psy} explains why the utilization ratio defined below is a pertinent mental workload metric. Here, the UDC is a human operator who has to service tasks from a queue. The types of services carried out by an operator include classification, supervision~\cite{Cummings2009Modeling-the-im,Supervisory} or assembly jobs within a production system~\cite{Ooijen2003The-effects-of-}.

\begin{definition} \label{def:UDC-UR} {\bf Utilization ratio:} The utilization ratio POCMC for a given positive averaging horizon $T$ is defined as:
\begin{equation}
W^{\text{\it \tiny R,$T$}}(k) =\frac{1}{T}\sum_{i=1}^{\min\{T,k\}} U(k-i), \quad k \in \mathbb{N}_+
\end{equation} where we adopt the convention that $\mathbb{W}^{\text{\it \tiny R,$T$}}=\mathbb{R}_+$, $W^{\text{\it \tiny R,$T$}}(0)=0$ and $U$ is a scalar non-negative process that governs the level of utilization. 
\end{definition}

In its simplest and most prevalent form, the set of utilization levels $\mathbb{U}$ would be $\{0,1\}$, and $U(k)=1$ and $U(k)=0$ would indicate whether the component is being used or not, respectively, at time $k$. The research on intelligent task management reported in~\cite{Savla2011A-Dynamical-Que}, also uses $\mathbb{U}=\{0,1\}$ for the following alternative mental workload metric quantifying utilization ratio with a forgetting factor. 

\begin{definition} \label{def:UDC-UR-with-forget} {\bf Utilization ratio with forgetting factor:} Given a forgetting factor $\alpha$ in $(0,1)$, the associated utilization ratio POCMC is defined as:
\begin{align}
S^{\text{\it \tiny  F},\alpha}(k+1) =& \ \alpha S^{\text{\it \tiny  F},\alpha}(k) + U(k), & k \in \mathbb{N} &, \  S^{\text{\it \tiny  F},\alpha}(0) = 0 \\ W^{\text{\it \tiny  F},\alpha}(k) = & \ (1-\alpha) S^{\text{\it \tiny  F},\alpha}(k), & k \in \mathbb{N} &
\end{align} where we adopt the convention that $\mathbb{W}^{\text{\it \tiny  F},\alpha}$ and $\mathbb{S}^{\text{\it \tiny  F},\alpha}$ are $\mathbb{R}_+$ and $U$ is a scalar non-negative process that governs the level of utilization. Notice that the performance process can be computed directly for $k$ greater than or equal to $1$ as $W^{\text{\it \tiny  F},\alpha}(k) =(1-\alpha)\sum_{i=0}^{k-1} \alpha^{k-i-1} U(i)$.
\end{definition}

 The authors of~\cite{Koch2009Channels-that-h} also adopt the utilization ratio with forgetting factor to model the dynamics of the temperature of the circuitry of the transmitter that broadcasts information across an additive noise link, subject to a transmission power process $U$. In this context, $\mathbb{U}$ is $\mathbb{R}_+$ and the UDC is the resulting communication channel whose performance is adversely affected by the thermal noise that intensifies with increasing temperature. Related work on allocation of energy harvested from a stochastic source for wireless transmission subject to constraints on temperature is reported in~\cite{Baknina2018aEnergy-harvesti}, while thermal effects were considered in~\cite{Forte2013Thermal-aware-s} in the context of distributed estimation.

 In contrast to the concept of mental workload, in the queueing literature the workload affecting the performance of the server quantifies the {\it effort} needed to complete the tasks apportioned to the server, but not yet completed. The following is a discrete-time approximation of the continuous-time model governing the workload process analyzed in~\cite{Baccelli2003Elements-of-que}. 

\begin{definition} {\bf Queueing Workload:} The following defines the queueing workload POCMC for a component acting as a server:
\begin{align}
S^{\text{\it \tiny  W}}(k+1) = &  S^{\text{\it \tiny  W}}(k) + U(k) - \check{Y}(k), \quad k \in \mathbb{N},\ S^{\text{\it \tiny W}}(0) = 0 \\
W^{\text{\it \tiny W}}(k) = & S^{\text{\it \tiny W}}(k)
\end{align} where $\mathbb{S}^{\text{\it \tiny  W}}$, $\mathbb{U}$ and $\mathbb{Y}$ are $\mathbb{N}$. Here, $U(k)$ and $\check{Y}(k)$ may represent the number of work quanta, or effort, associated with the incoming and completed tasks at time $k$, respectively. Notice that, given the structure in Fig.~\ref{fig:UDCModel}, $\check{Y}(k)$ must be determined as a function of $Y(k)$. We assume that $\check{Y}(k)$ must be zero when $S^{\text{\it \tiny  W}}(k)$ and $U(k)$ are zero, which requires a properly defined action kernel\footnote{In particular, this will require that $Y(k)$ carries enough information to determine when $U(k)$ is zero.}. 
\end{definition}  

\subsection{Energy harvesting models}
\label{subsec:EH_Models}

We also consider cases in which the state of the POCMC is governed not only by utilization but, unlike the models covered in Section~\ref{subsec:WorkloadUtilizationRatio}, also by extrinsic stochastic processes. Prime examples of these, include models of the so-called state of charge (SOC) quantifying the energy stored in a battery that is repeatedly recharged using energy harvested from unsteady sources. We propose the following model that is both a generalization and an adaptation to our discrete-time framework of ubiquitous models, such as those used in~\cite{Sudevalayam2011Energy-Harvesti,Kansal2007Power-Managemen}. Our model is general enough to capture the effects described in~\cite[Part~IV]{Priya2009Energy-Harvesti} for microbatteries that are often used in small devices powered by energy harvesting, including implanted medical devices~\cite{Lewandowski2009Feasibility-of-}.

\begin{definition}\label{def:EHmodel} {\bf Energy harvesting (EH) model} \\ Let $S^{\text{\it \tiny A}}:=\{ S^{\text{\it \tiny A}}(k) :  k \in \mathbb{N} \}$ be a given homogeneous Markovian process that quantifies not only the energy harvested over time but possibly also other stochastic phenomena that influence the operation of the battery and its recharging sub-systems. This process takes values in a subset $\mathbb{S}^{\text{\it \tiny A}}$ of a real coordinate space. A given map $\Delta^{\text{\it \tiny E}}:\mathbb{S}^{\text{\it \tiny SOC}} \times \mathbb{S}^{\text{\it \tiny A}} \times \mathbb{R}_+ \rightarrow \mathbb{R}_+$ governs the dynamics of the energy harvesting POCMC according to following recursion:
\begin{equation*} S^{\text{\it \tiny SOC}}(k+1) =  S^{\text{\it \tiny SOC}}(k) + \Delta^{\text{\it \tiny E}} \big (S^{\text{\it \tiny SOC}}(k), S^{\text{\it \tiny A}}(k),W^{\text{\it \tiny E}}(k) \big ), \quad k \in \mathbb{N}
\end{equation*} where we assume that $S^{\text{\it \tiny SOC}}(k)$, which quantifies the SOC at time $k$, is in $\mathbb{S}^{\text{\it \tiny SOC}}:=[ 0 , \bar{s}^{\text{\it \tiny SOC}} ]$ and $\bar{s}^{\text{\it \tiny SOC}}$ denotes the maximum SOC. The initial SOC is quantified by $S^{\text{\it \tiny SOC}}(0)$. Here, $W^{\text{\it \tiny E}}$ takes values in $\mathbb{R}_+$ and quantifies the energy effectively extracted from the battery for use by the action kernel. The map~$\Delta^{\text{\it \tiny E}}$ quantifies the net change in the battery charge resulting from the difference between the effect of $W^{\text{\it \tiny E}}(k)$ and the energy harvested. The state of the POCMC can be chosen as~$ S^{\text{\it \tiny E}} := \big \{ \big (S^{\text{\it \tiny SOC}}(k),S^{\text{\it \tiny A}}(k) \big ) \ : \ k \in \mathbb{N} \big \}$. The map $\mathcal{W}^{\text{\it \tiny E}}$ that determines $W^{\text{\it \tiny E}}$ in terms of $S^{\text{\it \tiny E}}$ and $U$ is described below in Remark~\ref{rem:mapWE}.
\end{definition}

The map $\Delta^{\text{\it \tiny E}}$ is characteristic of each battery and it must satisfy the following consistency conditions:
	\begin{align}
	\Delta^{\text{\it \tiny E}}(s^{\text{\it \tiny SOC}},s^{\text{\it \tiny A}},w^{\text{\it \tiny E}}) & \leq \ \bar{s}^{\text{\it \tiny SOC}} - s^{\text{\it \tiny SOC}}, \\ \Delta^{\text{\it \tiny E}}(s^{\text{\it \tiny SOC}},s^{\text{\it \tiny A}},w^{\text{\it \tiny E}}) & \geq \  -s^{\text{\it \tiny SOC}},
	\end{align} for all $s^{\text{\it \tiny SOC}}$, $s^{\text{\it \tiny A}}$ and $w^{\text{\it \tiny E}}$ in $\mathbb{S}^{\text{\it \tiny SOC}}$, $\mathbb{S}^{\text{\it \tiny A}}$ and $\mathbb{R}_+$, respectively. We also assume that $\Delta^{\text{\it \tiny E}}(s^{\text{\it \tiny SOC}},s^{\text{\it \tiny A}},w^{\text{\it \tiny E}})$ is continuous with respect to~$w^{\text{\it \tiny E}}$.

\begin{remark} \label{rem:mapWE} {\bf (Description of $\mathcal{W}^{\text{\it \tiny E}}$)} \\
As is known since the early work in~\cite{Shepherd1965Design-of-prima}, each battery type has a discharge curve that characterizes the voltage in terms of the SOC. Invariably, even for modern batteries~\cite{Chen2006Accurate-electr}, the voltage decreases as the SOC drops, which leads to the following constraints:
\begin{itemize}
	\item The maximum energy that can be delivered by the battery at any time $k$ is a decreasing function of the SOC, which we represent as $\mathcal{D}: \mathbb{S}^{\text{\it \tiny SOC}} \rightarrow \mathbb{R}_+$. More concretely, the constraint is given by:
	\begin{equation} \label{eq:ConstEHDischarge1}
	W^{\text{\it \tiny E}}(k) \leq \mathcal{D} \big ( S^{\text{\it \tiny SOC}}(k) \big ), \quad k \in \mathbb{N}
	\end{equation}
	\item There is a positive minimum SOC, denoted as $\underbar{s}^{\text{\it \tiny SOC}}$, below which the voltage is too low to power the component. This leads to the following constraint that must be satisfied for every $k$ in~$\mathbb{N}$:
	\end{itemize}
\begin{equation} \label{eq:ConstEHDischarge2}
	 \Big ( \Delta^{\text{\it \tiny E}} \big (S^{\text{\it \tiny SOC}}(k), S^{\text{\it \tiny A}}(k), W^{\text{\it \tiny E}}(k) \big ) + S^{\text{\it \tiny SOC}}(k) - \underbar{s}^{\text{\it \tiny SOC}} \Big ) W^{\text{\it \tiny E}}(k) \geq 0
	\end{equation}

Consequently, $U$, which represents the energy requested by the control policy, may differ from $W^{\text{\it \tiny E}}$. To be more specific, the energy extracted from the battery is determined in terms of $U$ and the state of the EH model via the map ${\mathcal{W}^E: \big ( s^{\text{\it \tiny E}}, u \big) \mapsto w^{\text{\it \tiny E}}}$ specified as follows:
\begin{multline*}
w^{\text{\it \tiny E}} = \max \Bigg \{\tilde{w}^{\text{\it \tiny E}} \geq 0 \ \Big | \ \tilde{w}^{\text{\it \tiny E}} \leq u, \ \tilde{w}^{\text{\it \tiny E}} \leq \mathcal{D} \big ( s^{\text{\it \tiny SOC}} \big ),  \\  \Big ( \Delta^{\text{\it \tiny E}} \big (s^{\text{\it \tiny SOC}}, s^{\text{\it \tiny A}},\tilde{w}^{\text{\it \tiny E}} \big ) + s^{\text{\it \tiny SOC}} - \underbar{s}^{\text{\it \tiny SOC}} \Big )  \tilde{w}^{\text{\it \tiny E}}  \geq 0 \Bigg \}
\end{multline*}
\end{remark}

The following is a simplified version of the EH model that is characterized by $\Delta^{\text{\it \tiny E}}$ with a linear range in which it quantifies the difference between the energy used and harvested, and it also implements a saturation that restricts the SOC to the interval~$\mathbb{S}^{\text{\it \tiny SOC}}$. 
\begin{definition} \label{def:LinSatEHmodel} {\bf (Linear-saturated EH model)} \\ The linear-saturated EH model is specified as follows for every $k$ in~$\mathbb{N}$:
\begin{subequations}
\begin{align} \label{eq:LinSatModel-a}
S^{\text{\it \tiny SOC}}(k+1) = & \min \{ S^{\text{\it \tiny SOC}}(k) + S^{\text{\it \tiny A}}(k) - W^{\text{\it \tiny E}}(k), \bar{s}^{\text{\it \tiny SOC}}\} \\
W^{\text{\it \tiny E}}(k) = & \begin{cases} U(k) & \text{if $U(k) \leq S^{\text{\it \tiny SOC}}(k) + S^{\text{\it \tiny A}}(k)$ } \\  S^{\text{\it \tiny SOC}}(k) + S^{\text{\it \tiny A}}(k) & \text{otherwise} \end{cases} 
\end{align}
\end{subequations} where $S^{\text{\it \tiny A}}(k)$, in this simplified model, represents the energy harvested at time $k$.
\end{definition}

The linear-saturated EH model does not capture the effects of the discharge curve and the changes that the mechanisms of charge and discharge go through as the SOC varies. The following is another simplified model in which the SOC takes values in a finite set and evolves as a controlled Markov chain (CMC).

\begin{definition} \label{def:SimpEHmodel} {\bf (Finite-state EH model)} \\
The SOC evolves according to a CMC whose state $S^{\text{\it \tiny E}}:=\{ S^{\text{\it \tiny E}}(k) : k \in \mathbb{N} \}$ takes values in ${\mathbb{S}^{\text{\it \tiny E}}:=\{0,1,\ldots,\bar{s}^{\text{\it \tiny E}} \}}$. The following map determines the probability transition map for $S^{\text{\it \tiny E}}$ in terms of $U$:
\begin{multline*}
\mathcal{S}^{\text{\it \tiny E}}(s^+|s,u,y) = \\ 
\begin{cases} \hat{\Gamma}_s(u) & \text{\it if $s^+=s+1$, $s < \bar{s}^{\text{\it \tiny E}}$} \\ \check{\Gamma}_s(u) & \text{\it if $s^+=s-1$, $s > 0$ } \\ 1-\hat{\Gamma}_s(u)-\check{\Gamma}_s(u) & \text{\it if $s^+=s$} \end{cases}, \\  \quad s^+ , s \in \mathbb{S}^{\text{\it \tiny E}}, \ u \in \mathbb{U}, \ y \in \mathbb{Y}
\end{multline*} where $\hat{\Gamma}_s: \mathbb{U} \rightarrow [0,1]$ and $\check{\Gamma}_s: \mathbb{U} \rightarrow [0,1]$ are given maps satisfying $\check{\Gamma}_0(u)=0$, $\hat{\Gamma}_{\bar{s}^{\text{\it \tiny E}}}(u)=0$ and $\hat{\Gamma}_s(u)+\check{\Gamma}_s(u) \leq 1$. Here, we assume that ${S}^{\text{\it \tiny E}}(k+1)$ is independent of $Y(k)$ when conditioned on ${S}^{\text{\it \tiny E}}(k)$ and $U(k)$ taking values $s$ and $u$, respectively. In this case, a given map $\mathcal{W}^{\text{\it \tiny E}}:\mathbb{S}^{\text{\it \tiny E}} \times \mathbb{U} \rightarrow \mathbb{R}_+$ determines the energy used by the action kernel as $\mathcal{W}^{\text{\it \tiny E}}: \big ( S^{\text{\it \tiny E}}(k),U(k) \big) \mapsto W^{\text{\it \tiny E}}(k)$.
\end{definition}


\section{Communication channels}
\label{sec:Channels}

In recent years, there has been a tremendous amount of research focused on energy harvesting wireless communication systems. For a comprehensive survey, we refer the reader to survey articles \cite{SudKul11, GunEtal14, UluEtal15}. As detailed above, an energy harvesting transmitter is a UDC with a state that indicates the amount of charge available for usage. In what follows, we provide a brief survey of two key areas: determining channel capacities and optimal scheduling policies in energy harvesting systems.

From an information-theoretic perspective, a key problem is identifying the capacity of an energy harvesting communication channel. The capacity of an additive white Gaussian noise (AWGN) channel with an energy harvesting transmitter was analyzed in  references \cite{OzeUlu12} and \cite{OzeUlu11}, for the infinite battery case and the no-battery case, respectively. Various upper and lower bounds on capacities have been studied in \cite{DonOzg12, JogAna14, TutEtal14, DonEtal15, OzeEtal15}. A general formula for the capacity of a point-to-point energy harvesting channel was established in \cite{ShaEtal16}. Reference \cite{ShaEtal16} also established a novel connection between the channel capacity and the optimal throughput (discussed below) for an energy harvesting transmitter. Beyond point-to-point channels, the capacity of energy harvesting multiple access channels (MACs) has been analyzed in \cite{OzeUlu12b, InaEtal18}, where a general capacity formula is derived, along with lower and upper bounds on capacity.

A significant amount of research has focused on the problem of scheduling for an energy harvesting transmitter. In this case, the transmitter has an energy queue as well as a data queue, and the goal is to transmit data to the recipients in the least amount of time, or equivalently transmit the maximum amount of data until a certain time. This problem has been studied in the offline setting, where the energy arrivals  are non-causally known, as well as the online setting where the transmitter has causal information about energy and data arrival \cite{GunEtal14, UluEtal15}. For the offline case, a variety of channel models have been investigated including point-to-point channels, broadcast channels \cite{OzeEtal11, OzeEtal12}, interference channels \cite{TutYen12}, and MAC channels \cite{YanUlu12}. The online case has also been studied for the point-to-point channel \cite{ShaOzg16}, the broadcast channel \cite{BakUlu16}, and the MAC channel \cite{BakUlu18}. We refer to \cite{BakUlu18} for a thorough list of references concerning online and offline scheduling in energy harvesting channels. In all the models described so far, the transmitter utilizes  energy for the sole purpose of transmission. Energy harvesting transmitters which expend energy on sensing, computing, communicating, and possessing imperfect batteries have been surveyed in \cite{UluEtal15}.

\subsection{Channels with evolving power constraints}

In addition to energy harvesting systems, we show that the UDC framework may also be used to analyze more general communication channels with time evolving power constraints. We describe these constraints below. The AWGN channel is one of the most popular channel models in information theory due to its relevance in practical applications. Evaluating the capacity of the AWGN channel under a variety of power constraints is a problem that has received much attention in the literature. The classical constraint studied by Shannon \cite{Sha48} involved an average power constraint of $P_{avg}$. In other words, if $(X(1), X(2), \dots, X(n))$ is the input to a channel, then it must satisfy
\begin{align*}
\frac{1}{n} \sum_{k=1}^n \norm{X(k)}^2 \leq P_{avg}.
\end{align*}
Shannon showed that the capacity of this channel is achieved using a random Gaussian codebook. In addition to the average power constraint, another practically relevant power constraint is the peak power constraint. A peak power constraint of $P_{peak}$ stipulates that every input $X$ to the channel should satisfy $\norm{X}^2 \leq P_{peak}$. Finding the capacity of this channel in the scalar case was first studied by Smith \cite{Smi71}. Smith showed that, although it is not possible to express the capacity in a closed-form expression, it may be calculated efficiently. The key observation in \cite{Smi71} was that the capacity is achieved by a \emph{discrete} input distribution that is supported on a finite number of atoms in $[- \sqrt{P_{peak}}, \sqrt{P_{peak}}]$. 

The flexible UDC framework allows us to model a variety of power-constrained channels. For example, the average power and peak power constraint may be restated as
\begin{align*}
\norm{X(k)}^2 &\leq \min\Big\{ k P_{avg} - \sum_{j=1}^{k-1} \norm{X(j)}^2, 
    \ P_{peak}\Big\}, \ k \in \mathbb{N}_+.
\end{align*}
It is evident that the power constraint on the $k$-th channel use depends not only on $P_{avg}$ and $P_{peak}$ but also on the symbols transmitted prior to time $k$. Therefore, this power constraint is \underline{utilization dependent}. We now provide a description of the UDC framework used for modeling a large class of power-constrained communication channels.

\begin{definition}[Evolving power constraints] \label{def:EPC}
An evolving power constraint is defined via a sequence of functions 
$\{ P_k : k \in \mathbb{N} \}$, where $P_k: \real_+^{k} \to \real_+$
such that $P_{k}(u_1, u_2, \dots, u_{k})$ determines the power constraint on the $(k+1)$-th transmission $X(k+1)$, where $u_l$ is the power of the $l$-th transmitted symbol, i.e., $u_l = \norm{X(l)}^2$, $l \geq 1$. 
\end{definition}

\begin{definition}[Evolving power constraint with state]\label{def: state}
An evolving power constraint with a state is characterized by three sequences of functions: (i) $\{f_k : k \in \mathbb{N}\}$ with $f_k : \real^{k} \to \real$, (ii) $\{\tilde f_k : k \in \mathbb{N}\}$, where $\tilde f_k: \real^2 \to \real$, and (iii) $\{p_k : k \in \mathbb{N}\}$ such that $p_k : \real \to \real$. These functions satisfy  the property that  $f_{k}(u_1, \dots, u_{k}) = \tilde f_{k}(f_{k-1} (u_1, \dots, u_{k-1}), u_{k})$, i.e., the value of function $f_k$  at time $k$ can be computed from that of function $f_{k-1}$ at time $k-1$ and $u_{k}$.

An evolving power constraint is said to have a state $f_{k}(u_1, \dots, u_{k})$ at time $k+1$ if the sequence of functions $P_k$ in 
Definition~\ref{def:EPC} may be written as $P_{k}(u_1, \dots, u_{k}) = p_{k}(f_{k}(u_1, \dots, u_{k}))$. Thus, the power constraint on the $k$-th symbol depends on the history of transmitted symbols up to time $k-1$ through the state at time $k$. 
\end{definition}

Evolving power constraints may be used to describe several power constraints studied in the literature. We provide a few examples below:

\begin{example}
Consider the standard average power constrained communication channel. Here, the constraint on the $k$-th symbol $X(k)$ is given by
\begin{align*}
\norm{X(k)}^2 \leq k P_{avg}- \sum_{j=1}^{k-1} \norm{X(j)}^2,
\end{align*}
for some fixed $P_{avg} > 0$. This power constraint may be characterized as an evolving power constraint with state, as follows: For $k \in \mathbb{N}$, define $f_k(x_1, \dots, x_{k}) = \sum_{j=1}^{k} \norm{x_j}^2.$ This definition satisfies the property that $f_{k+1}(x_1, \dots, x_{k+1})$ can be calculated using $f_k(x_1, \dots, x_k)$ and $x_{k+1}$; in particular, $\tilde f_{k+1}(u,v) = u+v$. The power constraint functions $P_k$ as $P_k(x_1, \dots, x_{k}) = (k+1)P_{avg} - f_k(x_1, \dots, x_{k})$ for every $k \in \mathbb{N}$. 
\end{example}

\begin{example} For an average power constraint of $P_{avg}$ coupled with a peak power constraint of $P_{peak}$, the only change from above is that $P_k(x_1, \dots, x_k) = \min(P_{peak}, (k+1)P_{avg} - f_k(x_1, \dots, x_{k}))$. 
\end{example}

\begin{example}
For a windowed-average power constraint over a window $T$, the state is the total energy expended over the last $T-1$ transmitted symbols. By allowing states to be vector valued in $\real^{T-1}$, this constraint is easily accommodated in Definition \ref{def: state}.
\end{example}

\begin{example}
A $(\sigma, \rho)$-power constraint is found to be relevant in energy harvesting applications as well as neuroscience. The $(\sigma, \rho)$-power constraint is defined as follows: Let $\sigma, \rho \geq 0$. A codeword $(x_1, x_2, \dots, x_n)$ is said to satisfy a \emph{$(\sigma, \rho)$-power constraint} if
\begin{equation}\label{eh: eq: def}
\sum_{j = k+1}^l x_j^2 \leq \sigma + (l-k)\rho ~\mbox{ for all } ~ 0 \leq k < l \leq n.
\end{equation}
The $(\sigma, \rho)$-power constraint essentially imposes a restriction on how bursty the transmit power can be, by constraining the total energy consumed over every interval to be approximately linear in the length of the interval. In energy harvesting communication systems, a $(\sigma, \rho)$-power constraint may be used to model a transmitter that harvests $\rho$ units of energy per unit time, and is equipped with a battery with capacity of $\sigma$ units, which is used to store unused energy for future transmissions. The $(\sigma, \rho)$-power constraints can be expressed equivalently by tracking a state parameter $\sigma_i$, that keeps track of the tightest constraint among the $k+1$ inequalities for $x_{k+1}$. The state function $\sigma_k$  evolves as follows: 
\begin{align*}
\sigma_{k}(x_1, \dots, x_{k}) = \min(\sigma, \sigma_{k-1}(x_1, \dots, x_{k-1}) + \rho - x_{k}^2).
\end{align*}
The power constraint function $P_k$ is defined as $P_k(x_1, \dots, x_k) = \sigma_k(x_1, \dots, x_k) + \rho$. Note that $\rho$ need not be constant over time, and such dependence or variability with respect to time is useful in modeling energy harvesting with arbitrary amounts of energy $\rho_k$ harvested at time $k$.
\end{example}

We define a UDC model that imposes evolving power constraints on an action kernel that is a communication channel.

\begin{definition}[POCMC component]
Let $\{(P_k, f_k, \tilde{f}_k) : k \in \mathbb{N} \}$ be an evolving power constraint with state. The POCMC component has state $S(k)$ that tracks the state of the power constraint at time $k$, i.e. $S(k) = f_{k-1}(u_1, \dots, u_{k-1})$. An input $U(k)$ indicates the desired power output for time $k$, i.e. the energy required to send symbol $X(k)$. The performance process $W(k)$ is equal to the power constraint imposed on $X(k)$. In other words,  $W(k) = P_{k-1}(u_1, \dots, u_{k-1}) = p_{k-1}(S(k))$  is a deterministic function of $S(k)$. The state at time $k+1$ satisfies $S(k+1) = \tilde f_{k}(S(k), u_{k})$.
\end{definition}

\begin{definition}[Action kernel]
The action kernel is a communication channel with input $X$ and output $Y$. At time $k$, the $k$-th symbol $X(k)$ is scheduled to be transmitted. The output of the channel $Y(k)$ depends on the input $X(k)$, the noise in the channel, as well as the performance process $W(k)$. Two natural cases to consider are:
\begin{enumerate}
\item
If $U(k) \leq W(k)$, then $X(k)$ is transmitted unaltered. Otherwise, $X(k)$ is rescaled to have power $W(k)$, i.e. $X'(k)$ is transmitted across the channel where $X'(k) = X(k) \times \sqrt{\frac{W(k)}{U(k)}}$.
\item
If $U(k) \leq W(k)$, then $X(k)$ is transmitted unaltered. Otherwise, there is no transmission, i.e. $X'(k)$ is transmitted across the channel where $X'(k) = 0$.
\end{enumerate}
\end{definition}

Note that since these are deterministic power constraints, the transmitter can calculate the power constraints on the $k$-th symbol in advance and ensure that $X(k)$ satisfies these power constraints. However, this is not possible when the power constraints are random. An example of random power constraint is the following. Consider a stochastic process 
$E := \{E(k) : k \in \mathbb{N} \}$, where 
$E(k)$'s are independent and identically distributed
random variables. We may now define the state as $(S(k), E(k))$, and the power constraint on the $k$-th symbol is computed via $p_k(S(k), E(k))$. The state evolution proceeds as $S(k+1) = \tilde f_k(S(k), u_k, E(k))$. This particular formulation is relevant to energy harvesting communication systems, discussed in Section \ref{subsec:EH_Models}. In the absence of any output $O(k)$, the transmitter has no way of modifying its $k$-th symbol to satisfy the power constraints. A variety of feedback settings are worth considering:
$O(k) = S(k)$ or $O(k) = (S(k), E(k))$. Additionally, the transmitter may also receive feedback from the receiver, i.e., $O(k)$ contains $Y(k)$. The capacity of energy harvesting systems with feedback has also been investigated in recent years, and it has been found that feedback increases capacity \cite{ShaEtal15}. In addition to random $E(k)$, yet another setting to consider is when the sequence of $E(k)$ is completely arbitrary, but is known to lie in some set. This is analogous to arbitrarily varying channels \cite{Bla60, CsiKor11} and may also be modeled using the UDC framework.

\subsection{Additional channel models}

Point-to-point communication under memoryless channels is widely studied and well-understood. However, there are several settings where this simple channel model is not sufficient. For example, we may consider channels that change with time, channels with memory, or channels that are simply unknown such as arbitrarily varying channels \cite{Bla60, CsiKor11}. We observe that many channel models of interest may be cast in the UDC framework by interpreting the channel as a  utilization dependent channel. Such a reinterpretation suggests generalizations that may be harder to arrive at directly. We briefly describe three examples of interest below: interference channels, finite-state Markov channels~\cite{GolVar96}, and channels with action dependent states~\cite{Wei10}.

\textbf{Inter-symbol interference:} Consider a discrete-time Gaussian inter-symbol interference channel as found in \cite{HirMas88}. For an input $X(k)$ at time $k$ in $\mathbb{N}_+$, the output $Y(k)$ of this channel is given by 
\begin{align}\label{eq: isi}
Y(k) = \sum_{i=0}^{k-1} \alpha^{i}X(k-i) + N(k),
\end{align}
where $N(k)$ is AWGN, and $\alpha \in (0,1)$. Define the state of the channel at time $k$ to be $S(k) := \sum_{i=1}^{k-1} \alpha^{i}X(k-i)$, and let $U(k) = X(k)$ and $W(k) = S(k)$. Then, we have 
\begin{align*}
S(k+1) &= \alpha (S(k) + U(k))  \ \mbox{ and } \\
Y(k) &= X(k) + W(k) + N(k).
\end{align*}
The first equation shows how the channel state $S(k)$ depends on the channel usage via previous inputs to the channel. The second equation highlights the effect of the performance process $W(k)$ on the action kernel, which is the communication channel from the input $X$ to the output $Y$.

\textbf{Channels with action-dependent states:} Weissman~\cite{Wei10} proposed a communication channel model with a state $S(k)$ at time $k$ which may be altered via an action $U(k)$.  This model includes communication settings where the channel state is known to the encoder and allows for novel settings, such as channels with a ``rewrite'' option. Effective communication over such a channel involves manipulating the channel state via actions and then coding for the channel. Let the state of the channel $S(k+1)$ at time $k+1$ depend on the input $U(k)$ according to the kernel $\mathcal{P}_{S(k+1) | S(k), U(k)}$. The action kernel is simply the communication channel where the output $Y(k)$ depends on the input $X(k)$ and the performance process $W(k)$, which is simply the channel state $S(k)$ in this case (i,e., $W(k) = S(k)$), according to the kernel $\mathcal{P}_{Y(k)|X(k), W(k)}$. This setup completely captures the model in~\cite{Wei10}. 

\textbf{Finite-state Markov channels:} The finite-state Markov channel, which was introduced in \cite{GolVar96}, consists of a channel with state $S(k)$ that evolves according to a Markov process, independently of channel inputs and outputs. This may be thought of as a utilization-independent channel where the state evolves according to a Markov kernel $\mathcal{P}_{S(k+1)|S(k)}$. A simple generalization allows us to model utilization-dependent channels using a CMC: the channel state $S(k)$ evolves as a CMC with an input $U(k) = X(k)$ in accordance with a transition kernel $\mathcal{P}_{S(k+1)|S(k), U(k)}$. The channel input is also $X(k)$, and the channel output is $Y(k)$ that is generated as per $\mathcal{P}_{Y(k)|X(k), W(k)}$ with $W(k) = S(k)$. There are a number of open problems concerning such channels: Is it possible to calculate channel capacity in closed-form? If the input $X = \{X(k): 
k \in \mathbb{N}\}$ is fixed to be Markovian, what is the maximum achievable channel capacity?


\section{Remote estimation across a packet-drop link powered by energy harvesting}
\label{sec:RemoteEstimation}

We begin this section by describing a UDC consisting of a packet-drop link powered by energy harvested and stored according to the models delineated in Section~\ref{subsec:EH_Models}. The apportionment of energy for transmission of information across the link over time is governed by a control process. We then proceed to discussing a few research themes in which the link is used in a remote estimation context.

\subsection{Packet-drop links powered by energy harvesting}
    \label{subsec:RE-EH}
    
At each time $k$, the link can either convey unerringly a symbol in $\mathbb{X}$ or a packet drop occurs. Implementation of the packet-drop link using wireless communication requires, for each $k$, that a codeword appropriately encoding $X(k)$ is placed for transmission across one or more physical channels. The transmission of a codeword will, in general, require multiple uses of each channel. A decoder at the receiver attempts to recover $X(k)$ and a packet drop occurs when it fails due to an outage caused by fading, interference or other detrimental effects. If $\mathbb{X}$ is infinite, such as when it is a real coordinate space, we assume that the codeword length is large enough to encode $X(k)$ with negligible quantization error.

\begin{definition}{\bf (EH packet-drop link)} 
\label{def:EHLink}
A packet-drop link comprises an action kernel whose output alphabet is $\mathbb{Y}:=\mathbb{X} \cup \{ \mathfrak{E} \}$, where $\mathfrak{E}$ indicates a packet drop. The input-output relationship is specified as follows:
\begin{equation}
Y(k) = \begin{cases} X(k) & \text{if $L(k) = 1$} \\ \mathfrak{E} & \text{if $L(k) = 0$}
\end{cases} 
\end{equation} where the link process $L$ indicates that there is a successful transmission when $L(k)=1$ and the packet is dropped otherwise. We assume that a map $\mathcal{L}:\mathbb{W}^{\text{\tiny \it E}} \rightarrow [0,1]$ characterizes $L$ probabilistically as follows:
\begin{equation} \label{eq:ProbPackDrop}
\mathcal{P}_{L(k) | S^{\text{\tiny \it E}}(k), U(K)} \big ( 0 | s^{\text{\tiny \it E}}, u  \big ) = \mathcal{L}(w^{\text{\tiny \it E}}), \quad k \in \mathbb{N}, \ s^{\text{\tiny \it E}} \in \mathbb{S}^{\text{\tiny \it E}}, \ u \in \mathbb{U}
\end{equation} which quantifies the probability of packet drop. Here, $W^{\text{\tiny \it E}}$ and $S^{\text{\tiny \it E}}$ are obtained from the EH model described in Definition~\ref{def:EHmodel} or a simplified version, such as the one specified in Definition~\ref{def:SimpEHmodel}. In addition, we assume that $L(k)$, $S^{\text{\tiny \it E}}(k+1)$ and $O(k)$ are conditionally independent given $S^{\text{\tiny \it E}}(k)$ and $U(k)$.
\end{definition}

In a wireless communication setting, an outage causing a packet drop occurs when fading, which is stochastic in general, attenuates the transmitted signal to a point that the received power is below a threshold needed for decoding~\cite{Goldsmith2005Wireless-commun,Tse2005Fundamentals-of}. The threshold depends on the codeword length, noise, interference characteristics~\cite{Ozarow1994Information-the} and it may also be stochastic. Here, we assume that fading and the transmission power are constant during the transmission of the codeword encoding $X(k)$. Moreover, $W^{\text{\tiny \it E}}(k)$ represents the total energy used attempting to transmit $X(k)$. Hence, $\mathcal{L}$, which quantifies the probability of outage given the transmission power as in~(\ref{eq:ProbPackDrop}), is a non-increasing function that can be determined on a case-by-case basis, such as in~\cite{Beaulieu2006A-closed-form-e}. 

\subsection{Design of remote estimation systems: problem definitions}

Henceforth, we prioritize the discussion of research on the design of remote estimation systems. Our choice is motivated not only by applications, such as monitoring of physical processes, but also by relevance for the design of control systems. 

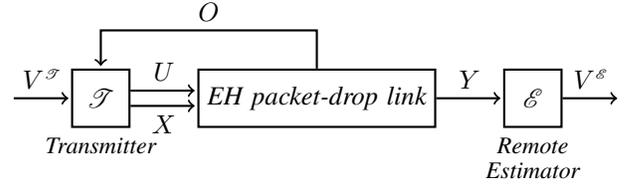
\begin{figure}
\centering

\begin{tikzpicture}[shorten >=1pt,on grid,auto]
    \tikzstyle{block}=[shape=rectangle,thick,draw,minimum size=0.3 in]

\node[block] (Transmitter) {$\mathscr{T}$};

\node[anchor=north] at (Transmitter.south) {\it \small Transmitter};

\node[right=0.35 in, block] at (Transmitter.east) (UDC) {\it EH packet-drop link };

\node[right=0.35 in, block] at (UDC.east) (Estimator) {$\mathscr{E}$};

\draw[->,thick] (UDC.east) -- (Estimator.west) node[midway,above] {$Y$};
\draw[->,thick] (Estimator.east) -- ([xshift=0.3 in]Estimator.east) node[midway,above] {$V^{\text{\tiny $\mathscr{E}$}}$};

\draw[->,thick] ([yshift=-0.04 in]Transmitter.east) -- ([yshift=-0.04 in]UDC.west) node[midway,below] {$X$};
\draw[->,thick] ([yshift=0.04 in]Transmitter.east) -- ([yshift=0.04 in]UDC.west) node[midway,above] {$U$} ;
\draw[->,thick] ([xshift=-0.3 in]Transmitter.west) -- (Transmitter.west) node[midway,above] {$V^{\text{\tiny $\mathscr{T}$}}$};

\draw[->,thick] (UDC.north) -- ([yshift=0.2 in]UDC.north) -- ([yshift=0.2 in]Transmitter.north) node[midway,above] {$O$}  -- (Transmitter.north);

\node[anchor=north] at (Estimator.south) (RemoteText) {\it \small Remote};
\node[anchor=north] at ([yshift=0.05 in]RemoteText.south) {\it \small Estimator};








\end{tikzpicture}
\caption{Basic overall structure of the remote estimation system considered.}
\label{fig:RemoteEstimationModel}
\end{figure}

We consider the configuration depicted in Fig.~\ref{fig:RemoteEstimationModel} in which an estimator $\mathscr{E}$ is a causal map that is possibly time-varying, and seeks to reconstruct a process $V$ based on information sent to it via a packet-drop link according to ${\mathscr{E}: y(1   \range  k) \mapsto v^{\text{\tiny $\mathscr{E}$}}(k)}$, for $k$ in $\mathbb{N}$. A transmitter is a causal map $\mathscr{T}$ that is possibly time-varying, and uses $V^{\text{\tiny $\mathscr{T}$}}$ and $O$ to produce $X$ and $U$ according to ${\mathscr{T}:(V^{\text{\tiny $\mathscr{T}$}}(1 \range k),O(1 \range k)) \mapsto (X(k),U(k))}$, for $k$ in $\mathbb{N}$. In most cases of interest $V^{\text{\tiny $\mathscr{T}$}}$ is either $V$ itself, or a causal function of $V$  possibly disrupted by additive or multiplicative noise. We refer to the pair $(\mathscr{T},\mathscr{E})$ in conjunction with the UDC that specifies the EH packet-drop link as a remote estimation system.


\begin{remark} \label{rem:EandTareInSync} {\bf Synchronization between $\mathscr{T}$ and $\mathscr{E}$} \\ 
Notice that one-step delayed feedback from the output of the link can be made available to $\mathscr{T}$ through $O$ by augmenting the state of the POCMC so as to include $Y(k-1)$. When such a feedback is present, a copy of the estimate $V^{\text{\tiny $\mathscr{E}$}}(k-1)$ can be replicated by $\mathscr{T}$ at time $k$. This {\it synchronization} often simplifies the joint design of $\mathscr{T}$ and $\mathscr{E}$ to meet stability or optimality conditions.
\end{remark}

We proceed with discussing the chronology of research on the design of remote estimation systems and control, with emphasis on the former.

\begin{problem}{\bf (Optimal remote estimation system design)} \label{prob:OptimalEst} \\ Let an EH packet-drop link, the joint probabilistic description of $V^{\text{\tiny $\mathscr{T}$}}(1 \range k)$ and $V(1 \range k)$ for all $k$ in $\mathbb{N}$ be given. For predetermined sets $\mathbb{T}$ and $\mathbb{E}$ of allowable transmitters and remote estimators, respectively, determine whether a pair $(\mathscr{T},\mathscr{E})$ exists that is optimal with respect to a given figure of merit ${\mathscr{J}:\mathbb{T} \times \mathbb{E} \rightarrow \mathbb{R}_+}$ that should assess the quality of $V^{\text{\tiny $\mathscr{E}$}}$ relative to $V$ and can include additional costs. If such a pair exists, determine one.
\end{problem}

Unless stated otherwise, we assume the following widely-used covariance-based cost structure:
\begin{multline}
\mathscr{J}^{\text{\tiny $(2q,K)$}}(\mathscr{T},\mathscr{E}) := \\ \frac{1}{K} \sum_{k=1}^K \EX \Bigg[ \Big( \big (V^{\text{\tiny $\mathscr{T}$}}(k)-V(k) \big )^T (V^{\text{\tiny $\mathscr{T}$}}(k)-V(k) \big ) \Big)^q \Bigg]
\end{multline} where $q$ is a positive integer and $K$ indicates the length of the optimization horizon. 

Stabilizability in the $m$-th moment sense, as defined below, is another relevant design objective.

\begin{problem}{\bf ($m$-th moment stabilizability)} \label{prob:StableEst} \\ Let an EH packet-drop link, the joint probabilistic description of $V^{\text{\tiny $\mathscr{T}$}}(1 \! \range \! k)$ and $V(1 \! \range \! k)$ for all $k$ in $\mathbb{N}$ be given. Consider that the $m$-th moment of $V$ is unbounded. Determine whether a pair $(\mathscr{T},\mathscr{E})$ exists for which the $m$-th moment of $V(k) - V^{\text{\tiny $\mathscr{E}$}}(k)$ is bounded for all $k$ in $\mathbb{N}$. If such a pair exists, determine one.
\end{problem}

Notice that the existence of a solution that is optimal for $\mathscr{J}^{\text{\tiny $(2q,K)$}}(\mathscr{T},\mathscr{E})$ in the limit when $K$ tends to infinity may imply, under certain conditions, $2q$-th moment stabilizability.

When either $\mathbb{T}$ or $\mathbb{E}$ is a singleton in Problems~\ref{prob:OptimalEst} or~\ref{prob:StableEst}, we say that the associated design problem is of the \underline{single-block} type, and we qualify it as \underline{two-block} otherwise.

\begin{remark}{\bf (Relevance of remote estimation for control systems)} There are at least two scenarios for which Problems~\ref{prob:OptimalEst} or~\ref{prob:StableEst} are relevant in the context of control systems. The first is when a packet-drop link connects the sensors that access the output of the plant to the controller. In this case, the transmitter is collocated with the sensors and the remote estimator is typically a component of the controller. The second setting is when the controller includes a transmitter to send its command signals to a remote estimator that is collocated with the actuator. A combination of both cases is also possible.
\end{remark}

\subsection{Uncontrolled transmission: optimal policies}
\label{subsec:RemoteEstimUncontrTrans}

As is surveyed in~\cite{Hespanha2007A-survey-of-rec}, the design of stabilizing and, whenever possible, optimal estimation and control systems whose components communicate via packet-drop links has been an active research topic for at least fifteen years. Early work assumed that the link process $L$ was an \underline{uncontrolled} time-homogeneous Markov chain. This assumption is realistic when the fading process, as indexed by $k$, is a real-valued time-homogeneous Harris chain and $\mathscr{T}$ does not have the authority to select the transmission power, which may be kept constant thanks to a dependable energy supply.

Henceforth, we limit our discussion to remote estimation systems in which $V$ and $V^{\text{\tiny $\mathscr{T}$}}$ are obtained as follows:
\begin{subequations}
\label{eq:V-LTIdynamics}
 \begin{align} \label{eq:V-LTIdynamics-a}
V(k+1)=&\mathbf{A} V(k) + N(k), \quad k \in \mathbb{N} \\ 
\label{eq:V-LTIdynamics-b}
V^{\text{\tiny $\mathscr{T}$}}(k) =& \mathbf{C} V(k) +  N^{\text{\tiny $\mathscr{T}$}}(k) , \quad k \in \mathbb{N}
\end{align}
\end{subequations} where $\mathbf{A}$ and $\mathbf{C}$ are real matrices of appropriate dimensions and the noise processes $N$ and $N^{\text{\tiny $\mathscr{T}$}}$ are independent and white with nonsingular covariance. In the context of control systems, an additional input term may be present in the right hand side of~(\ref{eq:V-LTIdynamics-a})~and~(\ref{eq:V-LTIdynamics-b}). 

At first, the effect of uncontrolled packet drops was modeled as multiplicative noise~\cite{Hadjiscostis2002Feedback-contro,Sinopoli2004Kalman-filering}, which makes the analysis of stability and second moment optimal design amenable to techniques inspired on Markovian jump linear system theory~\cite{Costa1993Stability-resul}. Typically the noise process would be Bernoulli, which would take value $0$ when a drop occurs. In a control systems setting, these multiplicative noises could affect the links carrying sensor measurements to the controller and control signals to the actuator. Most approaches focused on single-block design, which, depending on which links suffer packet drops, would be either a component at the sensors that processes measurements prior to transmission, a controller~\cite{Imer2006Optimal-control} or a remote estimator. As a consequence of the simplicity of the single-block framework, optimal policies and tight stabilizability conditions for state estimation and control can be obtained even when there is no link output feedback~\cite{Imer2010Optimal-estimat,Schenato2007Foundations-of-}, which can be viewed as a form of user datagram protocol~(UDP). 

The two-block remote estimation system formulated in~\cite{Xu2005Estimation-unde} was the first to consider the simultaneous design of $\mathscr{T}$ and $\mathscr{E}$. When~(\ref{eq:V-LTIdynamics}) is detectable~\cite{Rugh1996Linear-system-t,Hespanha2018Linear-systems-} the approach in~\cite{Xu2005Estimation-unde}, which is specified in continuous-time, can be immediately adapted to our discrete-time framework. In such a case, when $L$ is a Bernoulli process, the remote estimation system is $m$-th moment stabilizable if and only if the following condition holds:
\begin{equation} \label{eq:PacketdropUncontrolledStabilizability}
p^{\text{\it \tiny outage}} \rho(\mathbf{A})^{m} < 1 
\end{equation} where $p^{\text{\it \tiny outage}} : = \mathcal{P}_{L(k)}(0)$ is the probability of drop and $\rho(\mathbf{A})$ is the spectral radius of $\mathbf{A}$. As is shown in~\cite{Xu2005Estimation-unde}, a stabilizing solution is obtained by selecting $\mathscr{T}$ as a Kalman filter and $X$ as its state followed by a properly designed estimator $\mathscr{E}$. Subsequent work in~\cite{Gupta2007Optimal-LQG-con} showed that the scheme in~\cite{Xu2005Estimation-unde} is optimal with respect to a quadratic cost when $N$ and $N^{\text{\tiny $\mathscr{T}$}}$ are independent white Gaussian processes. Stabilizability in a control systems context was characterized in~\cite{Gupta2009Optimal-output-} using similar techniques for the case in which measurements are conveyed to the controller using two packet-drop links, with each having a distinct transmitter block. The setting in which a packet-drop link conveys command signals from the controller to the actuator was investigated in~\cite{Gupta2010On-stability-in}. 

Interestingly,~(\ref{eq:PacketdropUncontrolledStabilizability}) can be obtained as the limiting case~\cite{Minero2009Data-rate-theor} when $r$ tends to infinity of the condition in~\cite{Wong1997Systems-with-Fi,Martins2006Feedback-stabil,Sahai2006The-necessity-a} that characterizes stabilizability when a $r$-ary erasure channel\footnote{We refer the reader to the comprehensive overviews in~\cite{Matveev2009Estimation-and-,Hespanha2007A-survey-of-rec} that describe the most important classical results on control subject to data-rate constraints.} connects the transmitter to the remote estimator. 

\subsection{Controlled transmissions without packet drops}
\label{sec:ControlNoDrops}

We now consider the case in which transmissions may be controlled through $U$, while $V$ and $V^{\text{\tiny $\mathscr{T}$}}$ are modeled by~(\ref{eq:V-LTIdynamics}). When restricted to the remote estimation framework adopted here, in which $\mathscr{T}$ must designed to appropriately generate both $X$ and $U$, controlled transmissions were first studied in a stabilizability context in~\cite{Xu2005Estimation-unde}. 

In~\cite{Xu2005Estimation-unde}, $\mathscr{T}$ incorporates a Kalman filter that uses $V^{\text{\tiny $\mathscr{T}$}}$  to generate a local estimate $\hat{V}$ of $V$. In addition, it implements policies that use the magnitude of $\hat{V}-V^{\text{\tiny $\mathscr{E}$}}$ to determine the likelihood that a transmission is requested, which requires synchronization between $\mathscr{T}$ and $\mathscr{E}$ so that $V^{\text{\tiny $\mathscr{E}$}}$ can be reconstructed at the transmitter. Notice that, in the absence of packet drops, $\mathscr{T}$ and $\mathscr{E}$ can be synchronized without the need for feedback through $O$ since $Y$ can be causally computed at the transmitter based on $U$ and $X$. In this context, when $V$ is scalar and the noises $N$ and $N^{\text{\tiny $\mathscr{T}$}}$ are Gaussian, a policy that requests a transmission when the magnitude of $\hat{V}-V^{\text{\tiny $\mathscr{E}$}}$ exceeds a threshold was later shown\footnote{The techniques and results in~\cite{Martins2011Remote-state-es} are to a significant extent equivalent to the research reported in~\cite{Hajek2008Paging-and-regi} for paging and registration policies.} in~\cite{Martins2011Remote-state-es} to be optimal jointly with a Kalman-like estimator, with respect to a cost that linearly combines the expected squared estimation error and the time-averaged probability of transmission. As reported in~\cite{Nayyar2013Optimal-strateg}, threshold-type policies remain optimal when $V$ has dimension two or higher, provided that $A$ is a scaled orthogonal matrix. Although~\cite{Park2018Individually-op} shows that a jointly optimal transmitter and estimator pair exists for the aforementioned setting even when $A$ is any real-valued matrix, the question of whether there is a jointly optimal pair admitting threshold-type policies for transmission remains an open problem. Certainty equivalence properties for these estimators, which are relevant for the design of optimal controllers, are investigated in~\cite{Molin2013On-the-optimali}. Optimal strategies subject to restrictions on the total number of transmissions were determined in~\cite{Imer2006Optimal-control}. Results reported in~\cite{Wang2011Event-triggerin} show that threshold-based schemes can be adapted to guarantee stabilizability of a system formed by a network of plants and controllers connected by multiple packet-drop links. 

The framework in~\cite{Nayyar2013Optimal-strateg} was the first, in the context of remote estimation considered here, to allow for transmission policies that account for energy harvesting. Notably, it considers that $U$ is generated based not only on $V^{\text{\tiny $\mathscr{T}$}}$ but also on $S^{\text{\it \tiny SOC}}$, as determined by the linear-saturated EH model~(\ref{eq:LinSatModel-a}) for which the arrival process $S^{\text{\it \tiny A}}$ is assumed i.i.d. and $S^{\text{\it \tiny SOC}}$ is normalized so that each transmission at time $k$ requires $W^{\text{\it \tiny E}}(k)=1$. In this context, the following theorem, which follows from~\cite[Theorems~3 and~4]{Nayyar2013Optimal-strateg}, establishes an important structural result.

\begin{theorem} \label{thm:Martins1} When the noises in~(\ref{eq:V-LTIdynamics}) are zero-mean white Gaussian, and there are no packet drops, there are transmission and estimation policies with the structure in~(\ref{eq:KalmanLikeEst}) and~(\ref{eq:OptimalTransmNayyar}), respectively, that are jointly optimal for the scalar case.
\begin{equation}
\label{eq:KalmanLikeEst}
V^{\text{\it \tiny $\mathscr{E}$}}(k) = \begin{cases} A V^{\text{\it \tiny $\mathscr{E}$}}(k-1) & \text{if $Y(k)=\mathfrak{E}$} \\ X(k) & \text{if $Y(k)  \neq \mathfrak{E}$} \end{cases}, \quad k \in \mathbb{N}
\end{equation}

\begin{subequations}
\label{eq:OptimalTransmNayyar}
\begin{align} 
U(k+1) =& \begin{cases} 1 & \text{if $|\hat{V}(k+1)-AV^{\text{\it \tiny $\mathscr{E}$}}(k)| > $} \\ & \qquad \qquad \qquad \qquad \mathscr{G}(k,S^{\text{\it \tiny SOC}}(k+1)) \\ 0 & \text{otherwise} \end{cases} \\
X(k) =& \hat{V}(k) ,  \qquad k \in \mathbb{N} 
\end{align}
\end{subequations}
Here, $\mathscr{G}$ is a threshold that depends on time and the state of charge $S^{\text{\it \tiny SOC}}(k)$. The threshold determines when $U(k)$ is $1$, in which case a transmission setting $Y(k)$ to $X(k)$ is requested at time $k$. 
\end{theorem}
Methods to determine $\mathscr{G}$ are described in~\cite{Nayyar2013Optimal-strateg}. It is remarkable that a policy pair with the simple structure in~(\ref{eq:KalmanLikeEst}) and~(\ref{eq:OptimalTransmNayyar}) is jointly optimal, which also guarantees that it accomplishes the {\it best} trade-off between transmitting at time $k$ or saving energy to transmit later.

It is important to note that, barring the dependence of the thresholds for transmission on $S^{\text{\it \tiny SOC}}(k)$, (\ref{eq:KalmanLikeEst}) and (\ref{eq:OptimalTransmNayyar}) are akin to the optimal policies in~\cite{Martins2011Remote-state-es,Hajek2008Paging-and-regi}. As is explained in~\cite{Nayyar2013Optimal-strateg}, one way to obtain these results is to establish that there is a jointly optimal policy pair whose estimator has the structure~(\ref{eq:KalmanLikeEst}), after which the problem of finding a corresponding optimal transmission policy can be cast as an MDP~\cite{Puterman2005Markov-decision} whose state is finite dimensional because $\mathscr{T}$ and $\mathscr{E}$ are synchronized. Subsequently, well-known results can be invoked to prove that restricting transmission policies to be memoryless functions of the state of the MDP incurs no loss of optimality. Properties of the probability distributions of the noises, such as symmetry and unimodality, are used to show that there is no optimality loss when these  policies are further restricted to be of the form~(\ref{eq:OptimalTransmNayyar}).

\subsection{Controlled transmissions with packet drops and perfect feedback}
\label{subsec:RemEstContrTransFeed}

In this subsection, we discuss recent work for the framework that extends of that of Section~\ref{sec:ControlNoDrops} by allowing packet drops in the link that connects the transmitter to the remote estimator.

\begin{assumption}
\label{ass:Y-rec-from-O}
 Unless noted otherwise, here we assume that there is a causal map with which $Y(k)$ can be recovered unerringly from $O(1 \range k+1)$, for all $k$ in $\mathbb{N}$, which also implies that $\mathscr{T}$ and $\mathscr{E}$ can be synchronized.
\end{assumption}

\begin{assumption}
\label{ass:BatteryStateFrom-O} We also assume that $S^{\text{\it \tiny SOC}}(k)$ and $S^{\text{\it \tiny A}}(k)$ can be can recovered from $O(1 \range k)$, for all $k$ in $\mathbb{N}$.
\end{assumption}

We proceed to defining and subsequently discussing advantages and properties of a class of covariance-based transmission policies, which has been adopted in~\cite{Leong2018Transmission-sc,Li2017Power-control-o,Nourian2014Optimal-Energy-,Trimpe2014Event-based-sta}, to list a few.

A transmission policy is classified as covariance-based when the dependence of $U$ on $V^{\text{\it \tiny $\mathscr{T}$}}$ and $O$ can be recast in terms of a matrix-valued process $\mathbf{P}^{\text{\it \tiny $\mathscr{T}$}}$ that is determined from $Y$ as follows, for $k$ in $\mathbb{N}$:
\begin{equation}
\mathbf{P}^{\text{\it \tiny $\mathscr{T}$}}(k) :=\EX\Big[ \big ( V(k) - V^{\text{\it \tiny $\mathscr{E}$}}(k))^T ( V(k) - V^{\text{\it \tiny $\mathscr{E}$}}(k) \big) \ | \  Y(1 \range k) \Big]
\end{equation} where $\mathbf{P}^{\text{\it \tiny $\mathscr{T}$}}(0)$ is predetermined. 

There is a recursive time-update mechanism~\cite{Leong2018Transmission-sc} for $\mathbf{P}^{\text{\it \tiny $\mathscr{T}$}}$ that guarantees that it is an information state~\cite{Kumar2015Stochastic-Syst}, which, as we discuss below, may be used to recast the underlying optimization as an MDP, subject to the following set of policies.

\begin{definition} {\bf ( $\mathbb{T}^{\text{\it \tiny C}}$ - Memoryless~covariance-based~transmission~policy~set)}
\underline{We use $\mathbb{T}^{\text{\it \tiny C}}$} to denote the set of transmitters for which there is a map $\mathscr{T}^{\text{\it \tiny U}}$ determining $U$ according to ${\mathscr{T}^{\text{\it \tiny U}}:{\big ( P^{\text{\it \tiny $\mathscr{T}$}}(k-1),S^{\text{\it \tiny SOC}}(k),S^{\text{\it \tiny A}}(k) \big)} \mapsto U(k)}$.
\end{definition}

Now, consider the formulation in~\cite{Leong2018Transmission-sc}, in which, for each $k$, $\mathcal{T}$ selects $X(k)$ equal to $\hat{V}(k)$, and $U(k)$ is either zero (no transmission) or a pre-selected energy quantum, as opposed to allowing two or more energy levels. A transmitter that seeks to convey $X(k)=\hat{V}(k)$ to the estimator is often labeled {\it smart sensor} to distinguish it from the scheme in~\cite{Nourian2014Optimal-Energy-}, which attempts to forward the {\it unprocessed} measurements by setting $X(k)$ equal to $V^{\text{\it \tiny $\mathscr{T}$}}(k)$. 
 The following theorem establishing an important structural property for the estimator is a consequence of the analysis in~\cite[Section 2]{Leong2018Transmission-sc}.

\begin{theorem}\label{thm:SimpleEstimator} Assume that $N$ and $N^{\text{\it \tiny $\mathscr{T}$ }}$ are Gaussian and that we seek to solve Problem~\ref{prob:OptimalEst} subject to the additional constraint that the transmitter is in $\mathbb{T}^{\text{\it \tiny C}}$. In this context, restricting the set of estimators to be of the form (\ref{eq:KalmanLikeEst}) incurs no loss of optimality.
\end{theorem}

Consequently, Problem~\ref{prob:OptimalEst} becomes more tractable in exchange for the possible loss of optimality that results from restricting the transmitter structure to $\mathbb{T}^{\text{\it \tiny C}}$.
Notably, Theorem~\ref{thm:SimpleEstimator} allowed the authors of~\cite{Leong2018Transmission-sc} to show that there are coordinate-wise threshold transmission policies that are optimal among those in $\mathbb{T}^{\text{\it \tiny C}}$. In spite of these advantages, there is no known bound on the performance loss incurred by this method.

 The authors of~\cite{Li2017Power-control-o} investigated methods to determine optimal power selection policies when the probability of outage depends exponentially on the transmission power, which in their framework is allowed to vary among two or more levels. Notably, short of allowing for varying transmission power levels, the formulation of~\cite{Li2017Power-control-o} is analogous to the one in~\cite{Leong2018Transmission-sc}. Notwithstanding their similarities, the analysis in the former demonstrates why allowing the transmitter to select among multiple power levels complicates significantly the search and characterization of optimal policies. In order to contend with the complexity of the problem, work in~\cite{Li2017Power-control-o} includes useful approximations and tractable methods. The analysis and framework in~\cite{Knorn2017Optimal-energy-}, which also examines a control problem, provides suboptimal policies and numerical methods to address the case in which Assumptions~\ref{ass:Y-rec-from-O}~and~\ref{ass:BatteryStateFrom-O} are not satisfied.

Related work in~\cite{Ren2016Infinite-Horizo,Chakravorty2019Remote-estimati} extends~\cite{Lipsa2009Optimal-State-E} to the case in which information is sent from the sensor to the estimator across a noisy channel whose quality depends on an internal state and the transmission power. The research reported in these articles sought to solve Problem~1, in the absence of energy harvesting considerations, with respect to transmission policies that determine the transmission power $U$ and the transmitted signal $X$, and for a cost that combines the transmission power and the estimation error.

Tight necessary and sufficient conditions for the existence of a transmission policy that stabilizes the estimation error in the second-moment sense have been recently determined in~\cite{Lin2018Remote-state-es}. When adapted to our current formulation, assuming the finite-state EH model (Definition~\ref{def:SimpEHmodel}), the framework in~\cite{Lin2018Remote-state-es} would consider memoryless policies that use the state of charge $S^{\text{\it \tiny E}}(k)$ to decide, at each time $k$, whether a transmission should be attempted.
More precisely, the probability that a transmission is requested at time $k$ is a function of $S^{\text{\it \tiny E}}(k)$ represented as $\mathcal{L}^{\text{\it \tiny $\theta$}}:\mathbb{S}^{\text{\it \tiny E}} \rightarrow [0,1]$, which is denoted as $\theta$ in~\cite{Lin2018Remote-state-es}.
In order to employ the stabilizability conditions of~\cite{Lin2018Remote-state-es} in the current context, we assume that a map $\mathcal{L}^{\text{\it \tiny d}}:\mathbb{S}^{\text{\it \tiny E}} \rightarrow [0,1]$ is given which represents the probability of outage at time $k$ in terms of $S^{\text{\it \tiny E}}(k)$ when a transmission is requested at time $k$. Consequently, the map $\mathcal{L}^{\text{\it \tiny d}}$, which is represented with $d$ in~\cite{Lin2018Remote-state-es}, must account for the combined effect of the pre-selected policy that governs $U(K)$ in terms of $S^{\text{\it \tiny E}}(k)$ when a transmission is requested, $\mathcal{W}^{\text{\it \tiny E}}$ and $\mathcal{L}\big ( W^{\text{\it \tiny E}}(k) \big )$, which quantifies the outage probability according to~(\ref{eq:ProbPackDrop}).  Finally, we can restate~\cite[Theorem~3.1]{Lin2018Remote-state-es} in our context as follows. 

\begin{theorem} Given $\mathcal{L}^{\text{\it \tiny d}}$, there is a stabilizing transmission-request policy $\mathcal{L}^{\text{\it \tiny $\theta$}}$ if only if the following inequality holds:
\begin{equation}
\label{eq:stabilizabilityMichaelLin}
\lambda^{\text{\it \tiny S}} \rho(A)^2 \leq 1
\end{equation} where the nonnegative real constant $\lambda^{\text{\it \tiny S}}$ is a function of $\mathcal{L}^{\text{\it \tiny d}}$ and~$\mathcal{S}^{\text{\it \tiny E}}$. 
\end{theorem}

In addition, it is stated in~\cite[Theorem~3.1]{Lin2018Remote-state-es} that it suffices to consider deterministic transmission-request policies and according to~\cite[Theorem~3.2]{Lin2018Remote-state-es} the search can be further narrowed to threshold policies when $\mathcal{L}^{\text{\it \tiny $\theta$}}$ is non-increasing. 

Notice that (\ref{eq:stabilizabilityMichaelLin}) is a generalization of (\ref{eq:PacketdropUncontrolledStabilizability}) and the two conditions coincide when $\mathcal{L}^{\text{\it \tiny d}}$ is constant and equal to $p^{\text{\it \tiny outage}}$.


\section{Scheduling: Queueing Systems with Time-Varying
	Parameters, and Wireless Energy Transfer}
\label{sec:Queues}

The UDC model constitutes a natural framework 
for studying the scheduling problem in a wide
range of applications. In this section, we consider
two specific applications -- 
queueing systems with time-varying 
parameters and wireless systems with devices
powered by wireless energy
transfer -- and discuss how the UDC framework
can be used to study interesting and challenging
problems in these areas.

\subsection{Queueing systems with time 
    varying parameters}
	\label{subsec:QTV}
    
There exists a large volume of literature 
on queueing systems with time-varying 
parameters, dating back to the studies by 
Conway and Maxwell \cite{ConwayMaxwell62}, 
Jackson~\cite{Jackson63}, Yadin an Naor
\cite{Yadin63}, Gupta~\cite{Gupta67} 
and Harris~\cite{Harris67}, most of which 
focused on the state-dependent service 
rates. We refer a reader interested in a 
summary of earlier studies on queues
with state-dependent parameters to 
\cite{Dshalalow} and references therein. 
Although many, if not most, of these studies
can be carried out using the UDC framework
or its variant, 
here we focus on a more recent development in 
this area. 

The performance and management of humans
has been the subject
of many studies in the past, e.g., 
\cite{YerkesDodson, Edie54, 
Borghini14, Shunko17}. Recently, with 
rapid advances in information and sensor
technologies, human supervisory control (HSC) 
became an active research area
\cite{Supervisory, Sheridan97}. 
In HSC, human supervisors play a crucial role in 
the systems (e.g., supervisory control
and data acquisition (SCADA)) and at times are
required to process a large amount of information
in a short period in order to make critical 
decisions (e.g., a possible imminent nuclear 
meltdown due to a malfunction of cooling 
system), potentially
causing information overload. For this reason,  
there is a resurging interest in modeling
and understanding the performance of humans
under widely varying settings.
Although this is still an active research area, 
it is well documented that the performance
of humans depends on many factors, including
arousal and perceived workload
\cite{Edie54, Shunko17, Sheridan97, 
Asaro07, KcTerwiesch09}. 
For example, the well-known Yerkes-Dodson law
suggests that moderate levels of arousal are
beneficial, leading to the inverted-U model
\cite{YerkesDodson}. Moreover, the performance
of a human for varying types of tasks
(e.g., easy tasks vs. difficult tasks) changes
differently as a function of the level of 
arousal~\cite{Diamond07}. 

Savla and Frazzoli proposed a dynamical queue
approach to studying task management with human 
operators, using a differential system model
\cite{Savla2011A-Dynamical-Que}. The service
time of a task is equal to the product of 
(i) its workload and (ii) the value of service 
time function, which depends on the utilization 
level of the human operator at the time of
task assignment. The utilization level in 
their study is the continuous-time counterpart
of the utilization ratio with forgetting factor
defined in Section
\ref{subsec:WorkloadUtilizationRatio}. 
The service time function
determines the service time per unit workload
as a function of the utilization level. 

Their key results include the following:
\begin{enumerate}
\item[SF1.] Suppose that the service time function 
is convex and all tasks bring the same workload. 
Then, there is a maximally stabilizing  task release 
policy that applies a threshold to the
utilization level (Theorems III.1 and III.2
of \cite{Savla2011A-Dynamical-Que}): 
when the human operator
is idle at time $t$, the policy assigns a new
task if and only if its utilization level at time $t$
is less than or equal to some threshold. 

\item[SF2.] When tasks bring heterogeneous 
workloads, the maximum throughput that
can be achieved does not decrease compared
to the case with homogeneous workload of
tasks (Theorem IV.1 of 
\cite{Savla2011A-Dynamical-Que}). 
Thus, the heterogeneity of workload
does not diminish the maximum throughput
that can be achieved. 

\end{enumerate}

\subsubsection{UDC framework for 
	a Utilization-Dependent Server}

To the best of our knowledge, the work of
Savla and Frazzoli is the first to study the
task scheduling problem for servers whose
efficiency varies with an internal state
that summarizes its (recent) utilization.
Recently, Lin et al.~\cite{Lin-SQ, Lin-Workload} 
extended this work by adopting the UDC 
framework. Here we briefly summarize the work. 
A more detailed treatment is provided  
in Section~\ref{sec:UtilizationDependent}. 

The authors of \cite{Lin-SQ, Lin-Workload} 
introduced an internal state of the server, 
which is modeled using a finite-state CMC and 
approximates the utilization ratio
with forgetting factor. 
However, unlike in the model adopted in  
\cite{Savla2011A-Dynamical-Que}
where the service time depends only on the
utilization level at the time of task assignment, 
the service rate of the server is not assumed
fixed during the service time of a task; 
instead, it continues to change as the 
utilization level of the server
evolves over time in accordance with a CMC. 
Their key findings include the following:

\begin{enumerate}
\item[LF1.] Their study characterizes the maximum task arrival rate 
$\lambda^\star$, for which there exists a stabilizing task 
scheduling policy (see Theorem~\ref{thm:SQnecessary} in 
Section~\ref{sec:UtilizationDependent}). Moreover, the study
provides a computationally efficient method for computing 
$\lambda^\star$, 
by proving that this maximum task arrival rate is equal to 
the maximum throughput we can achieve using a threshold 
policy on the utilization level. 

\item[LF2.] In addition, the study shows that there exists 
a task scheduling policy with simple structure which can stabilize 
the system for any task arrival rate smaller than 
$\lambda^\star$. Specifically, a throughput-optimal 
task scheduling policy is a threshold policy
(see Theorem~\ref{thm:SQsufficient} in Section
\ref{sec:UtilizationDependent}); there exists
some threshold $\tau^*$ on the utilization level so that
the policy assigns a new task to the server if and only
if the queue is non-empty and the utilization level 
is below the threshold. 

\item[LF3.] Finally, for a fixed task arrival rate, the
study proposes
a method for finding a task scheduling policy that not only
stabilizes the system, but also aims to minimize the proportion
of time the server is requested to work on tasks, which 
the authors call the ``utilization rate" of the policy. 
In particular, they demonstrate
that the proposed method can produce a task scheduling
policy whose utilization rate can be made arbitrarily close 
to the minimum utilization rate achievable subject to the
stability constraint
(see Theorems~\ref{thm:MinWL} and \ref{thm:MinWLP}
in Section~\ref{sec:UtilizationDependent}).

\end{enumerate}

\subsection{Wireless energy transfer}

Wireless energy transfer (WET) has emerged as a potential solution  
to powering small devices that have low-capacity batteries  
or cannot be easily recharged, e.g., Internet-of-Things (IoTs) 
devices \cite{Bi16, Niyato17, Chu18, Lyu19}. 
Since the devices need to collect sufficient  
energy before they can transmit and the transmission rate is 
a function of transmit power, a transmitter has to decide 
(i) when to harvest energy and (ii) when to transmit and what
transmission rate it should use. 

For example, the studies reported  
in \cite{Ju14, Che15, Yang15} 
examined the problem of maximizing 
throughput in wireless networks in which communication 
devices are powered by hybrid access points 
(H-APs) via WET. 
In a related study, Shanet al.  \cite{Shan16} studied 
the problem of minimizing  the  total transmission delay  
or completion time of a given set of packets. 
In \cite{Rezaei19}, Rezaei et al. investigated
how to maximize the sum secrecy throughput 
among devices wirelessly powered by a base 
station and how to achieve max-min fair or 
proportionally fair secrecy throughput. 

In spite of recent efforts, there are many interesting
open problems. For example, what are the packet arrival 
rates at the devices powered by WET for which we can 
find a scheduling policy that stabilizes the queues at
the devices, subject to an energy expenditure 
rate or average power 
constraint, while capturing the battery 
discharge curve explained in Section
\ref{subsec:EH_Models}?
Another related question is: how do we  
minimize the energy expenditure rate at an H-AP
or a base station, 
while maintaining stable queues at the wireless 
devices and honoring (average) packet delay 
requirements? 

Here, we describe a model based 
on the UDC framework, which can be employed to
study these problems. For the simplicity 
of illustration, we focus on a scenario where 
two devices are powered by a single H-AP and
communicate on the uplink to the H-AP, and 
packets of fixed-size $K$ arrive at the devices 
according to  two independent Bernoulli processes 
with parameters $\boldsymbol{\lambda}
= (\lambda_0, \lambda_1)$.

\subsubsection{UDC framework for WET} 

Let $S^{soc} \Eqdef \{S^{soc}(k): k \in 
\mathbb{N} \}$, where $S^{soc}(k) = 
(S^{soc}_0(k), S^{soc}_1(k)) \in 
{\mathbb{S}^{soc}}^2$, 
$S^{soc}_i(k)$ represents 
the battery SOC of device 
$i$ at time $k$, and 
$\mathbb{S}^{soc}$ is the set of possible 
SOCs (defined in Section~\ref{subsec:EH_Models}).
Similarly, $Q = \{Q(k): k \in \mathbb{N} \}$, 
where $Q(k) = (Q_0(k), Q_1(k)) \in 
\mathbb{N}^2$ and $Q_i(k)$ is the number 
of backlogged bytes in the queue of device 
$i$. The packets of size $K$ can be 
segmented for transmission on the uplink from 
the devices to the H-AP. 

A scheduling policy is described by a mapping
$\theta: {\mathbb{S}^{soc}}^2 \times \mathbb{N}^2
\to \{H, 0, 1\} \times \mathbb{T}_{HAP} 
\times \mathbb{T}_{D}$, where 
$\mathbb{T}_{HAP}$ and $\mathbb{T}_D$
are the set of admissible transmit powers 
for the H-AP and the devices, respectively, 
with associated modulation and coding schemes.  
For a given pair $({\bf s}^{soc}, {\bf q})
\in {\mathbb{S}^{soc}}^2 \times \mathbb{N}^2$, 
$\theta_1({\bf s}^{soc}, {\bf q}) \in \{H, 0, 1\}$ 
indicates who will transmit, and 
$\theta_2({\bf s}^{soc}, {\bf q}) \in 
\mathbb{T}_{HAP}$
and $\theta_3({\bf s}^{soc}, {\bf q}) 
\in \mathbb{T}_{D}$ represent the selected 
transmit power: 
$\theta_1({\bf s}^{soc}, {\bf q}) = H$ 
means that the H-AP charges the devices 
using transmit power
$\theta_2({\bf s}^{soc}, {\bf q})$. On the other hand,  
$\theta_1({\bf s}^{soc}, {\bf q}) = i$ ($i = 0, 1$) 
informs device $i$ to transmit on the uplink
using transmit power $\theta_3({\bf s}^{soc}, 
{\bf q})$.  

When the H-AP transmits with power $p_H \in 
\mathbb{T}_{HAP}$, device $i$ receives energy 
at the rate $h_i \cdot p_H$ (per scheduling 
block duration), where $h_i$ is the transfer 
gain. The scenario in which the amount of 
transferred power is governed by a stochastic 
process can be handled analogously.   
In addition, when device $i$ is asked
to transmit at power $p_D$, it transmits
$g_i(p_D)$ bytes with transmit power
$\xi(p_D, s^{soc}_i)$, 
where $s^{soc}_i$ is the SOC of
device $i$ at the time. Here, 
$\xi(p_D, s^{soc}_i)$ models the power
that can be delivered by the battery as a
function of requested power and its SOC. 
The probability that the transmission from 
device $i$ will be unsuccessful
is determined by an error function 
$e_i(g_i(p_D), \xi(p_D, s^{soc}_i))$. 
  
Let $U = \{U(k): k \in \mathbb{N}\}$, where 
$U(k)$ is the scheduling decision at time $k$. 
The performance process is given by $W = 
\{W(k): k \in \mathbb{N} \}$, where 
$W(k) = 1 - e_i(g_i(p_D), \xi(p_D, s^{soc}_i))$
if $U(k) = (i, \cdot, p_D)$ with $i \in \{0, 1\}$, 
and $W(k) = 0$ if $U_1(k) = H$. The output
process of the UDC action kernel is $Y = 
\{Y(k): k \in \mathbb{N}\}$ with $Y(k) \sim$
Bernoulli($W(k)$).

The SOC and queue size of device $i$ ($i = 0, 1$) 
evolve as follows: 
\begin{eqnarray}
&& \hspace{-0.3in} S^{soc}_i(k+1)  
    \label{eq:WET1} \\
& = & \begin{cases}
	S^{soc}_i(k) - \xi(U_3(k), S^{soc}_i(k)) & \mbox{if } U_1(k) = i  \\
	S^{soc}_i(k) + h_i \cdot U_2(k) & \mbox{if } U_1(k) = H \\
	S^{soc}_i(k) & \mbox{if } U_1(k) = 1-i 
\end{cases}
    \nonumber
\end{eqnarray}
and
\begin{eqnarray*}
Q_i(k+1)
= Q_i(k) + K  \cdot B_i(k) - \check{Y}_i(k)
\end{eqnarray*}
where $B_i(k) \sim$ Bernoulli($\lambda_i$) 
and indicates whether or
not there is a new packet arrival at device $i$ at time
$k$, and $\check{Y}_i(k)$ is the number of 
bytes successfully transmitted from device $i$
to the H-AP at time $k$ and is given by 
 \begin{equation*}
\check{Y}_i(k)
= \begin{cases}
	g_i(U_3(k)) Y(k) & \mbox{if } U_1(k) = i, \\
	0 & \mbox{otherwise.}
\end{cases}
	\label{eq:WET2}
\end{equation*}
Here, we do not consider battery charge leakage, 
which may be an issue with cheap low-power 
remote sensing devices or IoT devices.
However, this can be handled by modifying 
equation \eqref{eq:WET1}.
For more details on battery charge leakage, we 
refer to Section~\ref{subsec:RealisticBattery}.

\begin{remark}
We note that the problem of designing a scheduling 
policy that can stabilize the queues at the wireless devices 
subject to constraints is similar to the problem investigated
by Lin et al. in \cite{Lin-SQ}. Moreover, finding a 
scheduling policy that minimizes the energy expenditure
rate at the H-AP is closely related to the problem of minimizing the
utilization rate with utilization-dependent server, which  
is studied in \cite{Lin-Workload} using the UDC framework. For these reasons, we expect that 
the UDC framework and tools used by Lin et al. 
in~\cite{Lin-SQ, Lin-Workload} will prove
to be instrumental to studying the wireless
systems with WET. 
We refer a reader interested in more detail 
to Section~\ref{sec:UtilizationDependent}. 
\end{remark}


As discussed in Sections~\ref{sec:Channels}
and \ref{sec:RemoteEstimation} (more 
specifically Section \ref{subsec:RE-EH}),
there is already extensive 
literature on wireless systems with energy
harvesting, 
including many studies that adopt an 
MDP formulation. 
For example, Kashef and Ephremides 
\cite{Kashef12} studied
the problem of maximizing the average number
of successfully transmitted packets over a 
time-varying channel by a wireless device
with EH. In their study, the channel state 
is modeled using a simple two-state 
time-homogeneous Markov chain (so-called 
Gilbert-Elliot model), and they derived
structural properties of optimal policies.
In \cite{Mao14}, Mao et al. examined a
similar problem of maximizing the average
amount of successfully transmitted data
over the course of a sensor life, 
where the sensor life is modeled using a 
geometric random variable. They proposed
an optimal energy allocation algorithm 
using the value iteration. 
Ahmed et al. \cite{Ahmed16} investigated
the problem of maximizing the data arrival
rate that can handled by the system, 
subject to delay constraints, where data
arrive with a fixed constant rate at each
time. They considered two different delay 
constraints 
-- average and statistical delay constraints
-- provided some structural properties of
optimal policies and an online
algorithm. 

Even though these studies are closely
related to those in WET in that, in both 
cases, device batteries need recharging 
through either EH or WET, and adopt a 
similar MDP formulation, 
there is also a key difference between them.
In the literature with EH, 
energy is often assumed to arrive
according to a deterministic process
or a random process, which is beyond the
control of decision makers. On the 
other hand, in the problem of WET, 
the energy delivery via WET is a part
of decision making process. Thus, the
scheduling policy must carry out a suitable
trade-off between recharging device
batteries and scheduling data/packet transmissions
in accordance with the objective of the
designer.

%

\section{Motivating Example: A Single Server with 
    Utilization-Dependent Service Rate} 
	\label{sec:UtilizationDependent}

In many cases of interest, the instantaneous 
performance or service rate of a server 
depends on their (recent) utilization. 
For example, as mentioned before, 
the efficiency of human operators 
is not constant and varies with
several factors, such as arousal and fatigue
\cite{YerkesDodson, KcTerwiesch09}. 
Thus, in many applications with human 
operators making critical decisions
(e.g., air traffic control and nuclear 
plant monitoring), 
it is important to take into account their
alertness and (instantaneous) efficiency, 
in order to improve the performance of
overall systems. The setting in which
the service rate of a server varies as 
a function of its recent utilization has 
been studied by Savla and Frazzoli in 
\cite{Savla2011A-Dynamical-Que}, using a
continuous-time differential system model, 
and their main findings are summarized 
in Section~\ref{subsec:QTV}. 

In this section, we use the recent work 
by Lin et al. \cite{Lin-SQ, Lin-Workload}
on a single server system to 
illustrate how the UDC framework was 
successfully used to 
facilitate the study of a similar system 
and help them extend the results in 
\cite{Savla2011A-Dynamical-Que}.
In particular, the UDC
framework is leveraged to expedite
the computation of
the maximum average task completion rate
or throughput that can be achieved by any 
stationary task scheduling policy (TSP)
and to design effective TSPs
\cite{Lin-SQ, Lin-Workload}.

To this end, they examine following two 
problems: in the first problem, they are
interested in designing a simple TSP that can 
stabilize the queue for any task arrival rate
for which there exists a stabilizing TSP. 
In the second problem, given a fixed
task arrival rate, they aim to 
devise a TSP that not only stabilizes the queue, 
but also minimizes the long-term proportion of 
time the server is requested to work on tasks.

\subsection{Setup: a UDC model}

Consider a queueing system with a single
server and an infinite first-in-first-out 
queue. The server is assumed non-preemptive; 
once it starts servicing a task, it
continues until the task is completed 
before taking on a new task or resting.

New tasks arrive according to a Bernoulli
process $B := \{ B(k):  k \in \mathbb{N}
\}$ with parameter $\lambda \in (0, 1]$:
$B(k) = 1$ if there is a new task arrival
at time $k$, the probability of which is
equal to $\lambda$, and $B(k) = 0$ 
otherwise. The Bernoulli random variables
$B(k), k \in \mathbb{N}$, are mutually
independent. Although
Bernoulli arrivals are assumed 
to simplify our discussion, 
more general arrival distributions 
(e.g., Poisson distributions) can be
handled with appropriate changes.

A new task that arrives at time $k$ becomes 
eligible for service, beginning at time 
$k+1$, and is put in the queue until the
server is ready to service it.
Each new task brings a (random) workload, 
which is the amount of service that the
server needs to provide in order to complete
the task.
They make a simplifying assumption that the 
workloads of tasks are modeled using
independent and identically distributed
exponential random variables. Thus, the
residual workload of a task that is 
currently being serviced has the same 
distribution as the workload of a queued 
task waiting for service. 

\begin{figure}[h]
\centering

      \tikzstyle{node} = [shape=circle, minimum size=0.5cm, circular drop shadow, text=black, very thick, draw=orange!55, top color=white,bottom color=orange!20, text width=0.3cm, align=center]
      \tikzstyle{building} = [draw,outer sep=7,inner sep=7,minimum size=57,line width=1, very thick, draw=black!55, top color=white,bottom color=black!20]

      \tikzstyle{start} = [  font={\huge\bfseries}, shape=circle, minimum size=0.5cm, circular drop shadow, text=black, very thick, draw=black!55, top color=white,bottom color=red!80, text width=0.5cm, align=center]
      \tikzstyle{goal} = [font={\huge\bfseries}, shape=circle, minimum size=0.5cm, circular drop shadow, text=black, very thick, draw=black!55, top color=white,bottom color=green!80, text width=0.5cm, align=center]

      \tikzstyle{default_edge} = [draw=orange!55, line width=2]
      \tikzstyle{greedy routing edge} = [->, draw=green!55, line width=5]
      \tikzstyle{face routing edge} = [->, draw=red!55, line width=5]

\begin{tikzpicture}[scale=0.8, every node/.style={transform shape}]
\draw[black, thick] (-1.05,-0.25) -- ++(2cm,0) -- ++(0,-0.8cm) -- ++(-2cm,0);
\draw[black, thick] (0.65,-0.25) -- +(0,-0.8cm);
\draw[black, thick] (0.35,-0.25) -- +(0,-0.8cm);
\draw[black, thick] (0.05,-0.25) -- +(0,-0.8cm);
\draw[black, thick] (-0.25,-0.25) -- +(0,-0.8cm);
\draw[<-,thick] (-1.2,-0.6) -- +(-20pt,0) node[left] {task arrivals};
\draw[->,thick] (2,-0.6) -- +(20pt,0) node[right] {completed tasks};
\node[align=center] at (0cm,-1.5cm) {Queue};


\draw[black, thick] (3.6,-3.25cm) circle [radius=1.5cm];
\draw[black, thick] (1.45,-0.65cm) circle [radius=0.5cm] node {\footnotesize $\mu^u(\cdot)$};

\node[align=center] at (1.5,1.3cm) {Scheduler};

\node(action) at (3.6,-2.6) {\small CMC $S^u$};
\draw[black, thick] (2.6,-3.5) circle (0.15);
\draw[->,thick] (2.7,-3.4) .. controls (2.9,-3.2) and (3.3,-3.2) .. (3.5,-3.4);
\draw[->,thick] (3.5,-3.6) .. controls (3.3,-3.8) and (2.9,-3.8) .. (2.7,-3.6);
\draw[black, thick] (3.6,-3.5) circle (0.15);
\draw[black, thick] (4.6,-3.5) circle (0.15);
\draw[->,thick] (3.7,-3.4) .. controls (3.9,-3.2) and (4.3,-3.2) .. (4.5,-3.4);
\draw[->,thick] (3.7,-3.6) .. controls (3.9,-3.8) and (4.3,-3.8) .. (4.5,-3.6);

\node[align=left] at (-1.6,-3.2) {
$Q(k)$ - queue length\\
$S^a(k)$ - server availability \\
$S^u(k)$ - utilization level};
\node at (-1.0,0.6) {\footnotesize$Q(k)$};
\node at (2.65,0.4) {\footnotesize$(S^u(k), S^a(k))$};
\node at (0.75,0.4) {\footnotesize$U(k)$};
\draw[thick]  (0.5,1.8) rectangle (2.5,0.8);

\draw[<-,thick] (1.2,-0.2) -- (1.2,0.8);
\draw[<-,thick] (1.6,0.8) -- (1.6,-0.2);
\draw[->,thick] (-0.6,1.2) -- (0.4,1.2);

\node [align=left] at (-2.5,-4.35) {
$U(k)$ - action};
\draw [dashed] (-3.5875,-2.3875) rectangle (0.3,-4);
\node [align=left] at (-2.4875,-2.1875) {System State};
\draw [black, thick](-0.6,-0.2) -- (-0.6,1.2);
\draw [dash pattern=on 2pt off 3pt on 4pt off 4pt](1.4,-1.2) -- (2.2,-4.4);
\draw [dash pattern=on 2pt off 3pt on 4pt off 4pt](2,-0.8) -- (4.8,-2);
\end{tikzpicture}
\caption{The UDC framework for the study 
by Lin et al. in \cite{Lin-SQ, Lin-Workload}.} 
\label{fig:Lin}
\end{figure}
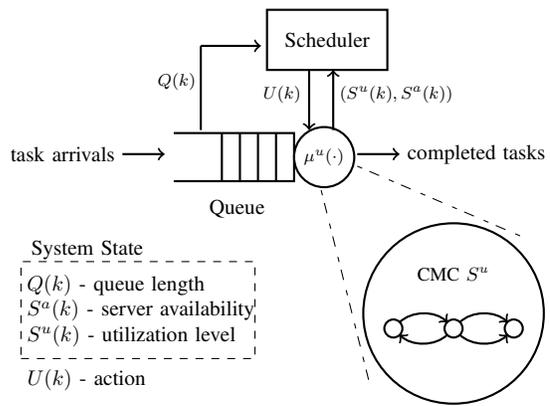

The control input at time $k$, namely
$U(k)$, is the scheduling decision 
chosen by the employed TSP:
$U(k) = 1$ means that the scheduler
requests the server to work on a 
task at time $k$, and $U(k) = 0$
otherwise, i.e., 
the server rests at time $k$. 
The control input process $\{U(k): k \in 
\mathbb{N}\}$ determines the evolution
of the server utilization over time, which in turn 
affects its instantaneous and long-term 
performance in the manner explained below. 

Let $S^u \Eqdef \{S^u(k): k \in \mathbb{N}\}$ 
be the process that tracks the utilization 
level of the server. For instance, $S^u(k)$ 
could represent the utilization ratio or the 
utilization ratio with forgetting factor 
$\alpha$ of the server (provided 
in Definitions~\ref{def:UDC-UR} and 
\ref{def:UDC-UR-with-forget} of Section
\ref{subsec:WorkloadUtilizationRatio}). 
In \cite{Lin-SQ, Lin-Workload}, 
the authors model $S^u$ using a
CMC taking values in a finite set
$\mathbb{S}^u := \{1, 2, \ldots, 
s_{\max}\}$. The dynamics of $S^u$
will be explained shortly.

The instantaneous performance or efficiency
of the server at time $k$ depends on its current
utilization level $S^u(k)$. This is modeled using
a service rate function $\mu^u: \mathbb{S}^u \to 
(0, 1)$, where $\mu^u(s^u)$ is 
the probability that
the server will complete a task within a 
unit time. The memoryless property
of an exponential distribution assumed for
task workloads implies that this
probability does not depend on the amount
of service that a task in service received
in the past. 

The transition probabilities of the utilization 
level $S^u(k) = s^u$ at time $k$
depend on (i) the current value of utilization, 
$s^u$, and (ii) the control input $U(k)$,
and are governed by the following mapping:
\begin{eqnarray*}
&& \hspace{-0.3in}
\mathcal{S}^u(\bar{s}^u | s^u, a)
:= \mathcal{P}(S^u(k+1) = \bar{s}^u |
	S^u(k) = s^u, U(k) = a) \\ 
& \hspace{-0.1in} = & \hspace{-0.1in}
\begin{cases}
	\eta^+_{s^u} 
		& \mbox{if } a = 1, 
			\bar{s}^u = \min(s_{\max}, s^u + 1) \\ 
	1 - \eta^+_{s^u}
		& \mbox{if } a = 1, 
			\bar{s}^u = s^u \\ 
	\eta^-_{s^u}
		& \mbox{if } a = 0, 
			\bar{s}^u = \max(1, s^u - 1) \\ 
	1 - \eta^-_{s^u}
		& \mbox{if } a = 0, 
			\bar{s}^u = s^u
\end{cases}
\end{eqnarray*} 
It is clear from the given transition probabilities
that if the server works on a task (resp. 
rests) at time $k$, the utilization 
level either remains at $s^u$ with probability 
$1 - \eta^+_{s^u}$ (resp. $1 - \eta^-_{s^u}$)
or goes up by one with 
probability $\eta^+_{s^u}$ if
$s^u < s_{\max}$ (resp. goes down by one
with probability $\eta^-_{s^u}$ if $s^u > 1$) 
with the convention 
$\eta_1^- = \eta_{s_{\max}}^+ = 0$.

The overall system dynamics are described by 
a CMC $S^{U} = \{(S^u(k), S^a(k), Q(k)): 
k \in \mathbb{N}\}$, where $Q(k)$ is the
number of backlogged tasks in the queue
at time $k$, and $S^a(k)$
indicates the availability of the server
to take on a new task. In other words, 
$S^a(k) = 1$ if the server
is available to service a new task at 
time $k$ (either after completing a
task or resting at time $k-1$), 
and $S^a(k) = 0$ otherwise. The CMC
$S^U$ takes values in the state space
of the system, which is given by 
$\mathbb{S}^U := \mathbb{S}^u \times 
\big( (\{0, 1\} \times \mathbb{N}) 
\setminus \{0, 0\} \big)$. 

In \cite{Lin-SQ, Lin-Workload}, they consider 
the following class of stationary
TSPs that map the current state of 
CMC, $S^U(k)$, to the probability of 
scheduling a task at each time $k$ in 
$\mathbb{N}$.

\begin{definition}
\label{defn:StationaryPolicy}
An admissible \underline{stationary 
randomized TSP} (SRTSP) is a mapping 
$\theta: \mathbb{S}^U \to [0, 1]$ such that 
(i) for all $(s^u, s^a, q) \in \mathbb{S}^U$,   
$\theta(s^u, s^a, q)$ is the probability
that the server is asked to work on a task 
when the CMC state is $(s^u, s^a, q)$, 
and (ii) $\theta(s^u, 0, q) = 1$ for all 
$q > 0$ and $s^u \in \mathbb{S}^u$.
\end{definition}

Obviously, the second requirement reflects
the assumption that the server is 
non-preemptive.
Also, under a fixed SRTSP $\theta$, the CMC 
$S^U$ is a discrete-time Markov chain with a 
countable state space. 
\\ \vspace{-0.05in}

\textbf{Threshold scheduling policies:} 
In practice, oftentimes a simple TSP is 
preferred as long as it does not cause
a significant degradation in performance. 
One class of simple
TSPs that are of interest is the set of
threshold TSPs:
fix a threshold $\tau \in \mathbb{S}^{u}_+
\Eqdef \{1, \ldots, s_{\max} + 1\}$. A 
threshold (task scheduling) policy with 
threshold $\tau$ is a deterministic
TSP given by a mapping
$\theta_\tau: \mathbb{S}^U
\to \{0, 1\}$, where  
\begin{eqnarray}  \label{eq:ThresholdPolicy} 
\theta_\tau(s^u, s^a, q)
\Eqdef \begin{cases}
	0 & \mbox{if (i) } s^u \geq \tau 
		\mbox{ and } s^a = 1  \\
	& \mbox{ or (ii) } q = 0, \\
	1 & \mbox{otherwise.}
	\end{cases}
\end{eqnarray}
Clearly, when the server is available to service 
a new task, the threshold policy $\theta_\tau$ 
assigns a new task if and only if the queue 
is non-empty and the utilization 
level is less than the threshold $\tau$.
Threshold policies are easy to implement in 
practice and require little information
for making scheduling decisions. 

When studying a queueing system, one of the most
important properties of interest is its stability;
an unstable system will lead to
poor performance in terms of the average
number of completed tasks per unit time 
or the (mean) sojourn times experienced by tasks.
For the study, they adopt the following notion
of stability.

\begin{definition} \label{defn:Stability}
For a fixed task arrival rate $\lambda > 0$, 
the CMC $S^U$ under a chosen 
SRTSP $\theta$, denoted by $S^U_\theta$, is 
said to be \underline{stable} if
\begin{enumerate}
\item there exists at least one recurrent 
communicating class of $S^{U}_\theta$;

\item all recurrent communicating classes are 
positive recurrent; and

\item the number of transient states is finite. 
\end{enumerate}
In addition, $\theta$ is said to 
\underline{stabilize} the CMC $S^U$ for 
the given task arrival rate $\lambda$. 
\end{definition} 
 
It is shown in \cite[Lemma 1]{Lin-SQ} that if 
$S^U_\theta$ is stable under some SRTSP $\theta$, 
there is a unique aperiodic, positive 
recurrent communicating class. As a result, 
we can find a unique stationary
distribution of $S^U_\theta$.

\subsection{Throughput-optimal task scheduling
	policies}
	\label{subsec:SQ-Stability}

As stated earlier, the stability of a system
is a fundamental property of a queueing 
system of interest. Therefore, a natural question 
that arises is: 
{\em how can we
design an effective TSP that stabilizes
the system whenever it is possible to do 
so, when the instantaneous performance of 
the server is affected by the very 
scheduling decisions it has made in the past?}

\subsubsection{Reduced process}

In order to find an answer to this question, Lin 
et al. investigated the problem of
designing a throughput-optimal TSP that 
stabilizes $S^{U}$ for any arrival rate 
$\lambda$ for which there exists a 
stabilizing TSP~\cite{Lin-SQ}. 
To this end, they first study a system in 
which there are infinitely many backlogged
tasks at the beginning, i.e., $Q(0) = \infty$. 
In this system, there is always 
a task waiting for service in the queue
when the server becomes available.

Consider the process 
$\tilde{S}^U = \big( (S^u(k), S^a(k) ): 
k \in \mathbb{N} \big)$ that describes
the utilization level and availability of
the server. It turns out that this 
reduced process $\tilde{S}^U$ 
plays a critical role in their studies.
 For example, 
when a threshold policy 
$\theta_\tau$ is adopted with some 
threshold $\tau > 1$, 
the resulting CMC $\tilde{S}^U$ 
can be modeled using a \underline{finite-state} 
Markov chain with a unique stationary distribution 
$\tilde{\pi}_\tau$ concentrated on the set 
\begin{eqnarray*}
\tilde{\mathbb{S}}_\tau
\Eqdef \{ (s^u, s^a) | s^u \in \{ \tau - 1, \ldots, 
s_{\max}\}, s^a \in \{0, 1\}  \}.
\end{eqnarray*}

As detailed in \cite{Lin-Workload} and summarized
below, this approach allows them to transform the
problem of designing a throughput-optimal 
TSP for $S^U$ to a more manageable MDP with a 
finite state space: the problem becomes one of
finding a TSP for the reduced process 
$\tilde{S}^U$ which maximizes the task completion 
rate of $\tilde{S}^U$. To tackle this problem,  
they leverage many well known results 
in the MDP literature, some of which are 
explained in Section~\ref{subsec:SQ-Discusion}. 

Define 
\begin{eqnarray}	\label{eq:lambda*}
\lambda^\star
\Eqdef \max_{\tau \in \mathbb{S}^{u+}} 
	\Big( \sum_{(s^u, s^a) \in 
		\tilde{\mathbb{S}}_\tau} 
		\tilde{\pi}_\tau(s^u, s^a) 
		\ \theta_{\tau}(s^u, s^a, 1) 
		\ \mu^u(s^u) \Big). 
\end{eqnarray}
Note that $\lambda^\star$ is the maximum 
average task completion rate 
among all threshold policies of the form in 
\eqref{eq:ThresholdPolicy} when the queue
is never empty. 
Let $\tau^\star$ be a maximizer of the 
right-hand side of \eqref{eq:lambda*}.

Clearly, threshold policies of the form 
in \eqref{eq:ThresholdPolicy} constitute
a small subset of the family of SRTSPs. 
Hence, without any additional assumptions, 
for instance, on the service rate function
$\mu^u$ as done in \cite{Savla2011A-Dynamical-Que}, 
one may suspect that we can find an SRTSP
that achieves a higher average task 
completion rate than $\lambda^\star$ and 
such an SRTSP will likely be able to 
stabilize the system for a task arrival 
rate larger than $\lambda^\star$. 
 
Somewhat surprisingly, this is not the case
and $\lambda^\star$, which 
can be computed efficiently by solving
the optimization in \eqref{eq:lambda*}
over a finite set, 
serves as an upper bound on the task arrival 
rate for which we can find a stabilizing 
SRTSP. This is formally stated by the 
following theorem. 

\begin{theorem}	\label{thm:SQnecessary}
Suppose that there exists a stabilizing 
SRTSP for some task arrival rate $\lambda 
> 0$. Then, $\lambda \leq \lambda^\star$. 
\end{theorem}

From the definition of $\lambda^\star$
in \eqref{eq:lambda*}, it is reasonable
to expect that when the task arrival rate 
$\lambda$ is smaller than $\lambda^\star$, 
we should be able to find a threshold 
policy that can stabilize the system, with 
$\theta_{\tau^\star}$ being a natural 
candidate. This is illustrated by the next 
theorem.

\begin{theorem} \label{thm:SQsufficient}
Suppose that the task arrival rate satisfies
$\lambda < \lambda^\star$. Then, the CMC
$S^U$ under the deterministic threshold
policy $\theta_{\tau^\star}$ is stable. 
\end{theorem}

An important implication of Theorems
\ref{thm:SQnecessary} and 
\ref{thm:SQsufficient} is that there 
is a throughput-optimal threshold 
TSP we can find efficiently.

\subsection{Utilization rate minimizing
	task scheduling policies}
	\label{subsec:SQ-Workload}

In some cases, in addition to keeping 
the system stable, it may be desirable
to minimize the proportion of time a server
is required to work. For instance, 
in the case of battery-powered wireless 
sensors sustained by renewable energy
or WET, 
we may wish to minimize the number of 
wireless transmissions over a long period,
while keeping the queue stable (assuming
that the measurements are not 
delay-sensitive). Similarly, when dealing
with human operators, it may be of 
interest to minimize the amount of
time a human operator is required to 
work on tasks. 

In a complementary study to \cite{Lin-SQ},
Lin et al. investigated the problem of designing a 
TSP that minimizes the utilization rate of 
the server, which is defined to be the 
long-term proportion of time the server works 
on tasks. This problem is considerably more 
challenging than the first problem as the
design requires maintaining system stability, 
while trying to minimize the utilization rate. 

For a given task arrival rate $\lambda \in 
(0, \lambda^*)$, let $\Theta(\lambda)$ 
be the set of stabilizing SRTSPs for $S^U$. 
Recall that, if a TSP $\theta$ stabilizes 
$S^U$, there is a unique stationary 
distribution of $S^U_\theta$, which we denote
by $\pi^\theta$. The corresponding 
utilization rate of the server in 
$S^U_\theta$ is equal to 
\begin{eqnarray*}
\mathcal{U}(\lambda; \theta)
\Eqdef \sum_{{\bf s} \in \mathbb{S}^U}
	\pi^\theta({\bf s}) \theta({\bf s}). 
\end{eqnarray*}

There are two questions of interest. 
First, given a task arrival
rate $\lambda$ in $(0, \lambda^\star)$, 
what is the minimum utilization rate we 
can achieve using a stabilizing TSP?
Denote this minimum utilization rate 
by  
\begin{eqnarray*}
\mathcal{U}^\star(\lambda)
\Eqdef \inf\{ \mathcal{U}(\lambda; \theta)
	\ | \ \theta \in \Theta(\lambda) \}.
\end{eqnarray*}
Second, how do we find a stabilizing
TSP that 
achieves a utilization rate arbitrarily
close to $\mathcal{U}^\star(\lambda)$?
In other words, given any constant 
$\delta > 0$, can we systematically 
find a stabilizing TSP 
$\theta^*$ such that
$\mathcal{U}(\lambda; \theta^*) \leq 
\mathcal{U}^\star(\lambda) + \delta$?

The authors of \cite{Lin-Workload} 
once again turned to the reduced process 
$\tilde{S}^U$ for answers: consider the 
problem of designing a TSP that makes a
decision 
on the basis of the utilization level and 
the availability of the server, assuming
that the queue is always non-empty (which
is true when $Q(0) = \infty$). 
Such a TSP is given by a 
mapping $\phi: \tilde{\mathbb{S}}^U 
\to [0, 1]$, where 
$\tilde{\mathbb{S}}^U 
\Eqdef \mathbb{S}^u \times \{0, 1\}$. 
As before, 
$\phi(s^u, s^a)$ represents the probability with
which the scheduler requests the server to work
on a task. We denote
the CMC $\tilde{S}^U$ under a TSP $\phi$
by $\tilde{S}^U_\phi$. 

They showed that, for every 
TSP $\phi$ satisfying $\phi(1, 1) > 0$, 
there is a unique stationary distribution 
of $\tilde{S}^U_\phi$. Denote this unique 
stationary distribution by 
$\tilde{\pi}^\phi$.\footnote{For a TSP $\phi$ 
with $\phi(1, 1) = 0$, there are two positive 
recurrent communicating classes of 
$\tilde{S}^U_\phi$ with one positive recurrent
communicating class being
$\{(1, 1)\}$.} Define $\Phi_+$ to be 
the set of TSPs $\phi$ satisfying 
$\phi(1, 1) > 0$. 
The utilization rate of $\phi$ in 
$\Phi_+$ is given by 
\begin{equation*}
\bar{\mathcal{U}}_\phi
= \sum_{(s^u, s^a) \in \tilde{\mathbb{S}}^U}
\big( \phi(s^u, s^a) \tilde{\pi}^\phi(s^u, s^a)
\big), 
\end{equation*}
and the average task completion rate
is equal to 
\begin{equation*}
\bar{\nu}_\phi
= \sum_{(s^u, s^a) \in \tilde{\mathbb{S}}^U}
\big( \phi(s^u, s^a) \tilde{\pi}^\phi(s^u, s^a)
\mu^u(s^u) \big).
\end{equation*}

Let $\Phi_+(\lambda)
\Eqdef \{ \phi \in \Phi_+ \ | \ 
	\bar{\nu}_\phi = \lambda \}$.

\begin{theorem}
	\label{thm:MinWL}
We have $\mathcal{U}^\star(\lambda)
= \bar{\mathcal{U}}_+(\lambda)$, 
where
\begin{equation*}
\bar{\mathcal{U}}_+(\lambda)
\Eqdef \inf\{ \bar{\mathcal{U}}_\phi
\ | \ \phi \in \Phi_+(\lambda) \}.
\end{equation*}
\end{theorem}

Although Theorem~\ref{thm:MinWL} reveals 
an interesting relation and a potential
means of computing $\mathcal{U}^\star(\lambda)$, 
it does not explicitly tell us how to 
find an SRTSP that can achieve a utilization
rate close to $\mathcal{U}^\star(\lambda)$. 
The following theorem sheds some light on 
this issue.

\begin{theorem}
	\label{thm:MinWLP}
Fix a task arrival rate $\lambda$ in 
$(0, \lambda^\star)$ and positive $\delta$. 
(i) We can find a 
pair $(\bar{\nu}^{(\lambda, \delta)}, 
\epsilon^{(\lambda, \delta)})$, 
where $\bar{\nu}^{(\lambda, \delta)} 
\in (\lambda, \lambda^\star)$ and 
$\epsilon^{(\lambda, \delta)}
\in (0, 1]$ such that 
\begin{equation*}
\Phi_{\epsilon}(\bar{\nu}^{(\lambda, \delta)})
\Eqdef \{ \phi \in 
	\Phi_+(\bar{\nu}^{(\lambda, \delta)})
	\ | \ \phi(1, 1) \geq 
	\epsilon^{(\lambda, \delta)}
	\}
\end{equation*}
is non-empty. (ii) Suppose $\phi^* \in 
\Phi_{\epsilon}(\bar{\nu}^{(\lambda, \delta)})$. 
Then, an SRTSP $\theta^*$ with 
\begin{eqnarray*}
\theta^*(s^u, s^a, q)
= \begin{cases}
	\phi^*(s^u, s^a) & \mbox{if } q > 0 \\
	0 & \mbox{otherwise}
\end{cases}
\end{eqnarray*} 
belongs to $\Theta(\lambda)$, and 
$\mathcal{U}(\lambda, \theta^*) 
\leq \mathcal{U}^\star(\lambda) + \delta$. 
\end{theorem}

Theorem~\ref{thm:MinWLP} indicates that
if we find a policy in 
$\Phi_{\epsilon}(\bar{\nu}^{(\lambda, \delta)})$, 
we can construct a stabilizing SRTSP 
whose utilization rate lies within $\delta$
of $\mathcal{U}^\star(\lambda)$. Hence, we
can get arbitrarily close to 
$\mathcal{U}^\star(\lambda)$ by reducing
$\delta$. Furthermore, they demonstrated in 
\cite{Lin-Workload} that the problem 
of finding a policy in
$\Phi_{\epsilon}(\bar{\nu}^{(\lambda, \delta)})$
can be formulated as a simple linear optimization
problem. Therefore, their study offers a systematic
way of constructing a suitable stabilizing
SRTSP whose utilization rate can be made
arbitrarily close to the minimum utilization 
rate.

\subsection{Discussion}
	\label{subsec:SQ-Discusion}

In this section, we first discuss the key
differences between the studies by Lin
et al. \cite{Lin-SQ, Lin-Workload} 
and that of Salva and Frazzoli
\cite{Savla2011A-Dynamical-Que}. Then, 
we outline how the UDC framework was vital
in obtaining the new results in 
\cite{Lin-SQ, Lin-Workload}. 

\subsubsection{Key differences between 
the studies by Savla and Frazzoli
\cite{Savla2011A-Dynamical-Que} and by Lin et al.
\cite{Lin-SQ, Lin-Workload}}

There are two key differences between these
two studies.
First, the authors of 
\cite{Savla2011A-Dynamical-Que} 
assume that the service time function 
is convex, which is analogous to the 
service rate function $\mu^u$ being unimodal
\cite{Lin-SQ, Lin-Workload}. 
But, Lin et al. do not impose any 
assumptions on the service rate function. 
In particular, $\mu^u$ is 
not assumed to be monotonic or unimodal 
and can be an arbitrary function taking values 
in (0, 1). Moreover, the model in 
\cite{Savla2011A-Dynamical-Que} 
assumes that the server efficiency during 
the service time of a task
depends only on the value of the 
service time function at the time the task was
assigned to the server and is
fixed during the service time of the task. In 
contrast, the model employed in \cite{Lin-SQ, 
Lin-Workload} allows the server 
efficiency captured by the service rate 
function to evolve while the server works on 
a task. 

Relaxing the assumptions introduced in 
\cite{Savla2011A-Dynamical-Que} and allowing 
a general service rate function is important 
to optimizing the performance of servers
with time-varying service rates, such 
as human operators, whose performance
is shown to be non-monotonic in arousal or
fatigue and its dependence on arousal 
varies with the difficulty level or
types of tasks, 
e.g., its performance tends to be 
monotonically increasing  
in arousal for easy tasks, whereas it is 
not the case for difficult tasks.

Second, in \cite{Savla2011A-Dynamical-Que}, 
a threshold policy is proved
to be maximally stabilizing only for 
the case with identical task workload. 
In the study by Lin et al., however, 
the workloads of tasks 
are modeled using i.i.d. random variables.
Although they assume an exponential distribution
for the workload to facilitate the analysis, 
similar results can be obtained with more 
general workload distributions with appropriate
changes to their model.

In addition to these key differences, there
is another important aspect of the results by
Lin et al. which
should be emphasized. Even though the results 
of \cite{Savla2011A-Dynamical-Que} are 
interesting, unfortunately they do not shed 
much light on a key practical question: 
{\em how do we identify a throughput-optimal 
policy in a computationally efficient manner?}
In contrast, the study by Lin et al. offers a 
systematic, 
computationally efficient method of finding a
throughput-optimal policy.

\subsubsection{Integral role of the 
UDC framework in \cite{Lin-SQ, 
Lin-Workload}}

It is noteworthy that the answers to several 
key questions by Lin et al. are obtained by 
relating the original problems to an MDP on a 
reduced, finite state space. This greatly simplifies 
their analysis and, more importantly, enables
them to leverage an extensive set of tools 
available for (constrained) MDPs. 

First, finding a throughput-optimal TSP in the 
UDC framework becomes straightforward;
an optimal threshold $\tau^\star$ can be 
computed as a solution to the optimization
problem in \eqref{eq:lambda*} by searching
through the \underline{finite} set 
$\mathbb{S}^{u+}$ with $s_{\max}+1$ 
elements. Note that, for each $\tau$
in $\mathbb{S}^u_+$, the stationary 
distribution $\tilde{\pi}_\tau$ has a
\underline{finite} support 
$\tilde{\mathbb{S}}_\tau$ and can be
computed efficiently as explained in 
\cite{Lin-SQ}. 

Second, the proofs of Theorems
\ref{thm:SQnecessary} through
\ref{thm:MinWLP} rely heavily 
on the tools available for
MDPs. For example, a well-known result
for constrained MDPs \cite[Theorem 
4.4]{AltmanMDP}
states that, for a constrained MDP, there
is an optimal policy that requires
at most $n_C$ randomizations, 
where $n_C$ is the number of constraints
in the corresponding constrained 
optimization problem.
Lin et al. formulate the problem of designing 
a throughput-optimal TSP as an unconstrained
MDP. Thus, they can infer from Theorem 4.4 of 
\cite{AltmanMDP} that there is an optimal 
deterministic TSP (which is not necessarily
a threshold TSP).\footnote{Although it is not
discussed in this article to keep our discussion
limited to single-queue scenarios, this observation 
plays an even more important role in an
extension of their work, where they
consider multiple types of tasks. The
problem of designing a stabilizing policy with, 
for example, two types of tasks 
can be viewed as one of maximizing the
throughput of one type of tasks 
\underline{subject to
a constraint} on the throughput
for the other type.}
In addition, the proof of Theorem
\ref{thm:MinWLP} requires computing an 
optimal policy for the MDP with a finite 
state space, which achieves
the minimum utilization rate among the 
policies in $\Phi_\epsilon(\nu)$ for 
some $\nu$ in $(\lambda, \lambda^\star)$. 
As this is an MDP on a finite state space, 
they are able to formulate an appropriate 
\underline{linear 
optimization problem}, the solution of 
which can be obtained efficiently and is 
used to construct an optimal policy that 
they seek. 

Third, we explain how the UDC framework
was central in bringing the synergy from
two different research fields in findings 
answers to their questions. Identifying 
near-optimal TSPs for the problem of 
minimizing the utilization rate 
demands expertise from optimization, 
queueing theory and stochastic control;  
finding such near-optimal TSPs calls for 
many structural results needed
to prove the main findings on stability. 
Their proofs borrow advanced tools from 
stochastic processes. These structural 
results are crucial to identifying 
optimal policies for MDPs on the reduced
state space via a linear optimization
problem, which are 
then used to construct near-optimal policies 
for the problem. Therefore, the UDP framework
is indispensable to carrying out the study 
at the intersection of these areas and 
obtaining the results. Finally, we 
point out that the problem studied
in \cite{Lin2018Remote-state-es} by
the same authors, was in part motivated by 
the UDC framework and, as a result, shares
a similar methodology.

\section{Conclusions and future directions}
\label{sec:ConcFutDir}
Our overview of the concepts, formulations, and methods utilized on the research themes expounded in Sections~\ref{sec:Channels}~-~\ref{sec:RemoteEstimation} evinces not only the similarities elicited by the presence of a UDC, but it also unveils a clear distinction among the objectives, techniques and assumptions adopted in each theme. This disconnection creates new research opportunities and challenges that would benefit from the fusion of the techniques and approaches that hitherto have been routinely employed by the information theory, wireless communication, operations research, networking and control theory communities. Notably, we concluded that the research challenges described in Sections~\ref{subsec:NoisyChannelsRemoteEstimation}-\ref{subsec:FeasibleRegionTradeoffPerformance} are currently not fully addressed, and constitute significant opportunities for future work that would also lead to methods for tackling problems specified by more realistic models and assumptions. Subsequently, in Sections~\ref{subsec:ManageSecurity}~-~\ref{subsec:RealisticBattery}, we proceed with suggesting additional future research directions that broach aspects of security and secrecy, effective methods to cope with systems comprising multiple UDCs, UDC in learning and more realistic battery models, respectively. 

\subsection{Noisy channels for remote estimation}
\label{subsec:NoisyChannelsRemoteEstimation}

Most work discussed in Section~\ref{sec:RemoteEstimation} presumes that, in the absence of an outage, an EH packet-drop link can convey a real vector unerringly from the sensor to the estimator when a transmission is requested. Future progress on new causal encoding and decoding schemes, possibly inspired on modifications of those discussed in Section~\ref{sec:Channels}, may lead to effective methods to tackle the unidealized case in which a noisy channel links the sensor to the estimator. Introducing channel encoding and decoding, and possibly lossy source compression, as was done in~\cite{Orhan2015Source-channel-} for an independent Gaussian source would expand the set of policies to include high and low fidelity solutions whose implementation may consume more or less energy~\cite{Grover2011Towards-a-commu}, respectively, in addition to that required for transmission. Obtaining methods for the design of such policies with stability and performance guarantees is, therefore, an important open challenge.

The case in which the UDC would depend not only on the energy available but also on the state of a physical system, such as the position and velocity vectors of a mobile agent, would be an interesting extension of this framework. In this setting, the UDC could be a communication channel between the agent and a base station whose outage likelihood would increase with distance for each transmission power level. The scenario in which the UDC would be a global positioning module (GPS) whose accuracy would depend on the location and power level, with higher fidelity consuming more power, would be an example relevant to autonomous navigation~\cite{Wheeler2018Relative-Naviga} of unmanned assets. In these cases, one needs to consider policies that not only allocate power for the UDC but also govern the control action that steers the agent. As is discussed in~\cite{Kreucher2005Sensor-manageme}, many active sensing problems could be formulated similarly once energy harvesting constraints are included.

\subsection{Queueing, remote estimation and age of information}
\label{subsec:QueueingRemoteEstimationAgeInfo}

According to the optimality principle used in~\cite{Gupta2007Optimal-LQG-con}, for the framework adopted in Section~\ref{subsec:RemoteEstimUncontrTrans}, if a sensor has access to the state $V$ or is able to compute the optimal state estimate $\hat{V}$~-~cases we refer to as {\it full-information} sensor or {\it smart sensor}, respectively~-~then it should always attempt to transmit the latest one to the remote estimator. Hence, given a choice, it is optimal to discard state estimates corresponding to failed transmission attempts in favor of the most recent one~-~a principle we term as {\it most-recent-only optimality}\footnote{See~\cite{Walsh2001Scheduling-of-N} for an overview of scheduling techniques for networked control systems, where most-recent-only scheduling is also discussed.}. In fact, this most-recent-only optimality principle for a full-information/smart sensor remains valid even in the controlled transmission setting described in Section~\ref{subsec:RemEstContrTransFeed}. Hence, these observations suggest that introducing a packet management layer, such as establishing a queue, prior to transmission is not necessary and may even be counterproductive when the sensor is full-state/smart. 

However, when using an existing transmission system one may be left with no option but to deal with a pre-existing first-in-first-out queue-based non-preemptive management system in which a packet leaves the queue only when it is successfully conveyed to the remote estimator. Notably, as is proved in~\cite{Sun2018Sampling-of-the,Sun2018Remote-Estimati} for the aforementioned scenario, for the case in which the source is a Wiener or Ornstein-Uhlenbeck process and the sensor is full-state, it is never optimal to submit a measurement for transmission when the queue is non-empty, and when a new measurement is inserted in the empty queue for transmission it must be the current state of the process, which can be viewed as a version of the {\it most-recent-only optimality} principle for the case when pre-emption is not allowed. Interestingly, the optimal rule proposed in~\cite{Sun2018Remote-Estimati} to determine whether to submit the latest measurement for transmission, subject to the queue being empty, follows an event-based threshold policy that is analogous to the one found to be optimal for the closely-related case analyzed in~\cite{Lipsa2009Optimal-State-E}\footnote{The techniques used in~\cite{Lipsa2009Optimal-State-E} are analogous to the ones adopted for the case without packet drops in~\cite{Martins2011Remote-state-es}.}. The fact that the most-recent-only optimality principle  may no longer hold when the sensor is neither full-state nor smart~\cite{Nourian2014Optimal-Energy-} raises the question of whether, if the sensor in the framework of~\cite{Sun2018Remote-Estimati} could transmit only noisy output measurements $V^{\text{\tiny $\mathscr{T}$}}$, there would be optimal policies for which a transmission would be scheduled even when the queue is non-empty. Furthermore, if the queue is served by a channel powered by energy harvested from stochastic sources then we are left with the currently unsolved problem of designing policies that determine not only when and which estimates or measurements should be placed in the queue for transmission but also allocate the energy used for each transmission attempt. A typical approach would be to characterize stabilizing policies first, perhaps within an appropriately parametrized class, followed by the characterization of structural properties that could facilitate the computation of optimal policies using tractable methods. The stability problem may require the integration of techniques such as the ones used in~\cite{Lin-SQ} and~\cite{Lin2018Remote-state-es}, which were discussed in Sections~\ref{sec:Queues} and~\ref{sec:RemoteEstimation} in the context of queueing and remote estimation, respectively. Devising methods to design optimal policies may involve fusing the techniques adopted in Section~\ref{subsec:RemEstContrTransFeed} and~\cite{Sun2018Sampling-of-the,Sun2018Remote-Estimati}, and possibly leveraging the fact that our UDC model is amenable to existing methodologies~\cite{Kumar2015Stochastic-Syst,Puterman2005Markov-decision} for POCMCs. Since the fidelity of the estimate constructed at the remote estimator depends on the recency of the information received by the remote estimator, both the stability and the optimization problems are related to recent work seeking to analyze and design data-transmission systems that effectively regulate the age of information~\cite{Sun2017Update-or-wait:,Arafa2018Age-minimal-tra}. In fact, it has been suggested in~\cite{Kam2018Towards-an-effe,Ornee2019Sampling-for-re} that the remote estimation and age of information problems are inextricably tied.


\subsection{Feasible region and trade-off among performance 
metrics}
\label{subsec:FeasibleRegionTradeoffPerformance}

Most existing studies in which queue length, utilization 
or workload affects the performance of servers, 
including those mentioned
in Sections~\ref{sec:Queues}
and \ref{sec:UtilizationDependent}, 
examine the effects on a
single aspect of server performance, oftentimes their 
service rates being the choice. In another example, 
the study by Chatterjee et al.~\cite{Chatterjee17} takes 
into account the service quality (which is modeled as 
channel condition in their study) as a function of queue 
length and examines the information-theoretic capacity 
of such systems. 

In many cases of interest, however, including human 
supervisors~\cite{Sheridan97}, several performance
aspects, including service rate and service quality
(e.g., reliability or frequency of mistakes or poor 
decisions), can be affected at the same time by work 
history via server state. Moreover, the 
requirements (e.g., service rate vs. reliability)
in different applications are likely to vary 
considerably based on the types of tasks that
need to be processed. 

From this viewpoint, it is important to develop a 
comprehensive theory for these systems, including 
their fundamental limits.
Regrettably, to the best of our knowledge, little
is known about the feasible region of multiple performance 
metrics which can be achieved simultaneously 
and how to design suitable policies for carrying out a 
desired trade-off among various performance metrics in the
feasible region, in particular on the Pareto frontier.

\subsection{Secure remote estimation powered by energy-harvesting}
\label{subsec:ManageSecurity}

Preventing, or at the very least mitigating the effect of, attacks on the channels connecting the sensors to every component relying on remotely constructed state estimates is critical to ensure the safe operation~\cite{Pajic2017Design-and-impl} of networked cyber-physical systems. While clever encoding and decoding schemes~\cite{Fawzi2014Secure-estimati}, some of which may be implemented efficiently using event-based algorithms, may thwart or curb the effect~\cite{Shoukry2016Event-triggered} of certain types of attacks, a relentless surreptitious Man-in-the-Middle (MitM) attack injecting false data~\cite{Mo2010False-data-inje} may significantly degrade the performance of any remote estimation system. The case-study in~\cite{Jovanov2018Secure-state-es} illustrates that by employing message authentication codes (MACs), even if infrequently, may afford performance guarantees against MitM attacks. It further demonstrates that although MAC are known to substantially increase communication overhead, which is particularly critical when using bandwidth-limited networks such as the ones found in automobiles, its parsimonious use may suffice for practical purposes. A promising new research avenue is to investigate estimation-oriented encoding and decoding schemes and MAC scheduling policies that would jointly provide stability and performance guarantees, or would even be jointly optimal with respect to a given estimation error metric, in the presence of MitM attacks. Realistic problem formulations, in which information transmission is powered by an energy harvesting module, would have to account for the additional energy required for the transmission of MAC. A new type of EH link\footnote{Possibly based on a modification of Definition~\ref{def:EHLink}.} that would account not only for packet-drop events, but also MitM attacks whose likelihood and severity would depend on the power employed in each transmission for the inclusion of MAC could be a useful abstraction to design and evaluate the performance of such systems.  The open problems discussed here would also be relevant for distributed function calculation~\cite{Sundaram2011Distributed-fun} in the cases in which information would be wirelessly disseminated among the agents via such security-threatened EH links. 

Finally, it would be important to investigate all of these problems in light of other security threats~\cite{Urbina2016Survey-and-new-}, including denial-of-service attacks~\cite{Cetinkaya2019An-overview-on-}. 



\subsection{Systems with multiple UDCs and
	development of a new theory for large
	systems with many UDCs}
\label{subsec:MultipleUDC}

In many situations of practical interest, there 
are a set of servers working on tasks (e.g., 
emergency rooms at hospitals). 
Furthermore, the availability
of servers may be affected by some exogenous processes
(e.g., schedules of doctors and nurses at hospitals). 
For example, data centers comprise a large number
of server racks that are connected by high-speed
networks and are sometimes subject to power 
constraints. Also, because the reliability 
of hardware components, such as CPUs, GPUs and 
memory modules, degrades when 
the temperature exceeds some threshold, they
need to be cooled for stable operation, for instance,
via direct-to-chip liquid cooling. Moreover, 
because new server racks are added over time 
to meet increasing demands and old or failed racks 
are replaced at different times, the computational 
capabilities offered by various computational resources, 
which are designed for different types of tasks
(e.g., CPUs vs. GPUs), can vary significantly.

Another class of problems well suited for the
UDC framework with multiple UDCs, 
which is also related to those 
in Sections~\ref{subsec:QueueingRemoteEstimationAgeInfo}
and \ref{subsec:ManageSecurity},
is information collection from multiple sources 
over time. These sources 
may be distributed sensors in wireless sensor
networks (WSNs), which are powered by renewable 
energy or WET, 
or ``friends" in social networks who
prefer not to be bothered constantly for 
the latest information.
One can view the ``usefulness" of the 
information collected from each sensor
or friend as the reward. Such usefulness
of information from a sensor or a friend will 
likely be stochastic. However, 
there are certain factors that would affect the 
usefulness of the information. These include
(i) the accuracy or quality of the sensors 
or the importance of the friends in social
networks (which are often measured using their 
``centralities" in social networks
\cite{EasleyKleinberg}) and 
(ii) the age-of-information from each sensor
or friend introduced in Section
\ref{subsec:QueueingRemoteEstimationAgeInfo}
as well as the frequency of information 
requests.

Unfortunately, the quality of sensors and
the importance of friends may not be known 
in advance. In addition, 
in many practical scenarios, we may be
able to poll or collect information from 
only a limited number of 
sensors or friends at any given time 
and only so often. In WSNs, for instance, 
the number of available 
channels or timeslots in a frame may constrain
the number of measurements we can collect
at each time and, when the sensors are 
powered by renewable energy or WET, they may not
be able to report measurements even when
they are polled, as their availability for
reporting measurements will be governed by 
a stochastic process. 

Despite their prevalence, not much is known about 
their fundamental performance limits and 
efficient resource management in such systems
with many UDCs. This calls for the development
of a new theory, especially for offering a 
guideline for effective resource allocations
in large systems with many UDCs. 
Only recently research has demonstrated
the benefits of task-aware scheduling at data
centers (e.g.~\cite{Dogar14}). Consequently, there 
is a rich set of open problems in related
domains. 
For instance, when heterogeneous servers are 
designed/optimized for different types of tasks
and their efficiency is utilization-dependent, 
how should we schedule arriving tasks so that 
both the (mean) sojourn time of the tasks and 
the utilization rates of the servers are minimized 
while maintaining the stability of the queues? 
These are some of questions, the answer
to which can have significant impact on many 
areas, including crucial applications 
involving HSC (e.g., 
air traffic control and nuclear power
plant monitoring). A useful approach for
studying these problems, especially 
when some of the parameters are unknown, 
is the restless multi-armed bandit model, which 
has been previously applied to stochastic
scheduling~\cite{Mahajan2008Multi-armed-ban, 
Whittle88}.



\subsection{UDC in learning}

In recent years, information-theoretic techniques have emerged as effective tools to study optimization procedures in machine learning problems \cite{RusZou16, XuRag17, PenEtal18, NegEtal19, SteEtal20}. Iterative and noisy optimization procedures such as the stochastic gradient descent or stochastic gradient Langevin dynamics have been hypothesized to be efficient due to their inherent noisy nature \cite{PenEtal18, GeEtal15, JinEtal17, AbbSan20}.  It is not hard to describe such iterative optimization techniques using the UDC framework, but such a reformulation leads to some novel adversarial models of learning that may be of interest. For example, the controlled Markov process may be current hypothesis $S(k)$ that is being optimized in an iterative fashion. The action kernel may be used to generate the additive update $Y(k)$ that modifies $S(k)$ to $S(k+1) = S(k) + Y(k)$. In general, this update is simply a derivative (or a noisy version thereof). In an adversarial setting, these updates may be influenced by an adversary via an input $X(k)$ to the action kernel. Such a scenario generalizes data-poisoning attacks studied in \cite{BigEtal12, MeiZhu15, KohLia17} or gradient-based attacks \cite{BlaEtal17, CheEtal19}.

\subsection{More realistic battery models for energy harvesting: leakage and nonlinearities}
\label{subsec:RealisticBattery}

Although, as we discussed in Section~\ref{subsec:EH_Models}, the batteries used in energy harvesting modules have a rather complex behavior, the existing work discussed throughout this article adopts either the linear-saturated or the finite-state approximations. These simplified models do not capture a host of issues that could possibly require new methods and abstractions. This is illustrated by the following two features that could be captured by our general model of Definition~\ref{def:EHmodel}:
\paragraph{\it Leakage} The chemistry of every battery and the operation of its auxiliary circuitry will cause charge to leak, even when it is not supplying power. Hence, the charge that is stored in a battery may be partially lost unless it is used quickly or the leakage is offset by harvesting.
This is a relevant problem for low-power remote sensing devices that operate over long periods of time.

\paragraph{\it Nonlinearities} In Section~\ref{subsec:EH_Models} we mentioned the fact that, due to the discharge curve, in general there is a state of charge threshold below which the voltage of the battery does not suffice to power the other components. Consequently, if the voltage is near the required minimum then leakage effects may drain the state of charge below the aforesaid threshold, after which enough energy must be harvested before the battery can function again. The fact that the state of charge also governs the portion of the energy harvested that is effectively stored constitutes another important nonlinearity. Notably, as the state of charge nears its maximum and minimum the ability of the battery to store energy varies considerably.

Although these detrimental battery features have been taken into account in the design of transmitters that seek to maximize wireless transmission rate subject to energy harvesting~\cite{Devillers2012A-General-Frame,Tutuncouglu2015Optimum-Policie}, they have yet to be addressed in the context of remote estimation systems.


\section*{Acknowledgement}
The authors would like to thank Sennur Ulukus (UMD), Yasser Shoukry (UMD) and Vijay Gupta (UND) for helpful discussions and suggestions. 

\bibliographystyle{ieeetr}
\bibliography{MartinsRefs-v3.bib,VarunRefs.bib}

\begin{thebibliography}{100}

\bibitem{Ulukus2015Energy-harvesti}
S.~Ulukus, A.~Yener, E.~Erkip, O.~Simeone, M.~Zorzi, P.~Grover, and K.~Huang,
  ``Energy harvesting wireless communications: a review of recent advances,''
  {\em IEEE Journal on Selected Areas in Communications}, vol.~33,
  pp.~360--381, March 2015.

\bibitem{Kim2012Cyber-physical-}
K.-D. Kim and P.~R. Kumar, ``Cyber-physical systems: a perspective at the
  centenial,'' {\em Proceedings of the IEEE}, pp.~1287--1308, May 2012.

\bibitem{Baheti2011The-impact-of-c}
R.~Baheti and H.~Gill, {\em The impact of control technology},
  ch.~Cyber-physical systems, pp.~161--166.
\newblock IEEE Control Systems Society, 2011.

\bibitem{Staal2004Stress-cognitio}
M.~A. Staal, ``Stress, cognition and human performance: A literature review and
  conceptual framework,'' Tech. Rep. NASA/TM-2004-212824, NASA, August 2004.

\bibitem{Leong2018Optimal-control}
A.~S. Leong, D.~E. Quevedo, and S.~Dey, {\em Optimal control of energy
  resources for state estimation over wirless channels}.
\newblock Briefs in eletrical and computer engineering, Springer, 2018.

\bibitem{Hespanha2007A-survey-of-rec}
J.~P. Hespanha, P.~Naghshtabrizi, and Y.~Xu, ``A survey of recent results in
  networked control systems,'' {\em Proceedings of the IEEE}, vol.~95, no.~1,
  pp.~138--162, 2007.

\bibitem{Teigen94}
K.~H. Teigen, ``{Yerkes-Dodson: a law for all seasons},'' {\em Theory \&
  Psychology}, vol.~4, no.~4, pp.~525--547, 1994.

\bibitem{Wickens2000Engineering-psy}
C.~D. Wickens and J.~G. Hollands, {\em Engineering psychology and human
  performance}.
\newblock Prentice Hall, third edition~ed., 2000.

\bibitem{Cummings2009Modeling-the-im}
M.~L. Cummings and C.~E. Nehme, ``Modeling the impact of workload in network
  centric supervisory control settings,'' in {\em Proceedings of the 2nd Annual
  Sustaining Performance Under Stress Symposium}, February 2009.

\bibitem{Supervisory}
J.~R. Peters, V.~Srivastava, G.~S. Taylor, A.~Surana, M.~P. Eckstein, and
  F.~Bullo, ``Human supervisory control of robotic teams: Integrating cognitive
  modeling with engineering design,'' {\em IEEE Control Systems Magazine},
  vol.~35, pp.~57--80, December 2015.

\bibitem{Ooijen2003The-effects-of-}
H.~P.~G. van Ooijen and J.~W.~M. Bertrand, ``The effects of a simple arrival
  rate control policy on throughput and work-in-process in production systems
  with workload dependent processing rates,'' {\em International Journal of
  Production Economics}, vol.~85, pp.~61--68, 2003.

\bibitem{Savla2011A-Dynamical-Que}
K.~Savla and E.~Frazzoli, ``A dynamical queue approach to intelligent task
  management for human operators,'' {\em Proceedings of the IEEE}, vol.~100,
  pp.~672--686, March 2012.

\bibitem{Koch2009Channels-that-h}
T.~Koch, A.~Lapidoth, and P.~P. Sotiriadis, ``Channels that heat up,'' {\em
  IEEE Transactions on Information Theory}, vol.~55, pp.~3594--3612, August
  2009.

\bibitem{Baknina2018aEnergy-harvesti}
A.~Baknina, O.~Ozel, and S.~Ulukus, ``Energy harvesting communications under
  explicit and implicit temperature constraints,'' {\em IEEE Transactions on
  Wireless Communcations}, vol.~17, pp.~6680--6692, October 2018.

\bibitem{Forte2013Thermal-aware-s}
D.~Forte and A.~Srivastava, ``Thermal-aware sensor scheduling for distributed
  estimation,'' {\em ACM Transactions on Sensor Networks}, vol.~9,
  pp.~53:1--53:31, July 2013.

\bibitem{Baccelli2003Elements-of-que}
F.~Baccelli and P.~Br\'{e}maud, {\em Elements of queueing theory}.
\newblock Springer, second edition~ed., 2003.

\bibitem{Sudevalayam2011Energy-Harvesti}
S.~Sudevalayam and P.~Kulkarni, ``Energy harvesting sensor nodes: Survey and
  implications,'' {\em IEEE Communications Surveys and Tutorials}, vol.~13,
  no.~3, pp.~443--461, 2011.

\bibitem{Kansal2007Power-Managemen}
A.~Kansal, J.~Hsu, S.~Zahedi, and M.~B. Srivastava, ``Power management in
  energy harvesting sensor networks,'' {\em ACM Transactions on Embeded
  Computing Systems}, vol.~6, September 2007.

\bibitem{Priya2009Energy-Harvesti}
S.~Priya and D.~J. Inman, eds., {\em Energy Harvesting Technologies}.
\newblock Springer, 2009.

\bibitem{Lewandowski2009Feasibility-of-}
B.~E. Lewandowski, K.~L. Kilgore, and K.~J. Gustafson, ``Feasibility of an
  implantable, stimulated muscle-powered piezoelectric generator as a power
  source for implanted medical devices,'' in {\em Energy Harvesting
  Technologies} (S.~Priya and D.~J. Inman, eds.), ch.~15, pp.~389--404,
  Springer, 2009.

\bibitem{Shepherd1965Design-of-prima}
C.~M. Shepherd, ``Design of primary and secondary cells ii. an equation
  describing battery discharge.,'' {\em J. Electrochem. Soc.}, vol.~112, no.~7,
  pp.~657--664, 1965.

\bibitem{Chen2006Accurate-electr}
M.~Chen and G.~A. Ric\'{o}n-Mora, ``Accurate electrical battery model capable
  of predicting runtime and i-v performance,'' {\em IEEE Transactions on Energy
  Conversion}, vol.~21, pp.~504--511, June 2006.

\bibitem{SudKul11}
S.~Sudevalayam and P.~Kulkarni, ``Energy harvesting sensor nodes: Survey and
  implications,'' {\em IEEE Communications Surveys \& Tutorials}, vol.~13,
  no.~3, pp.~443--461, 2011.

\bibitem{GunEtal14}
D.~Gunduz, K.~Stamatiou, N.~Michelusi, and M.~Zorzi, ``Designing intelligent
  energy harvesting communication systems,'' {\em IEEE communications
  magazine}, vol.~52, no.~1, pp.~210--216, 2014.

\bibitem{UluEtal15}
S.~Ulukus, A.~Yener, E.~Erkip, O.~Simeone, M.~Zorzi, P.~Grover, and K.~Huang,
  ``Energy harvesting wireless communications: A review of recent advances,''
  {\em IEEE Journal on Selected Areas in Communications}, vol.~33, no.~3,
  pp.~360--381, 2015.

\bibitem{OzeUlu12}
O.~Ozel and S.~Ulukus, ``Achieving {AWGN} capacity under stochastic energy
  harvesting,'' {\em IEEE Transactions on Information Theory}, vol.~58, no.~10,
  pp.~6471--6483, 2012.

\bibitem{OzeUlu11}
O.~Ozel and S.~Ulukus, ``{AWGN} channel under time-varying amplitude
  constraints with causal information at the transmitter,'' in {\em 2011
  Conference Record of the Forty Fifth Asilomar Conference on Signals, Systems
  and Computers (ASILOMAR)}, pp.~373--377, IEEE, 2011.

\bibitem{DonOzg12}
Y.~Dong and A.~{\"O}zg{\"u}r, ``Approximate capacity of energy harvesting
  communication with finite battery,'' in {\em 2014 IEEE International
  Symposium on Information Theory}, pp.~801--805, IEEE, 2014.

\bibitem{JogAna14}
V.~Jog and V.~Anantharam, ``An energy harvesting {AWGN} channel with a finite
  battery,'' in {\em 2014 IEEE International Symposium on Information Theory},
  pp.~806--810, IEEE, 2014.

\bibitem{TutEtal14}
K.~Tutuncuoglu, O.~Ozel, A.~Yener, and S.~Ulukus, ``Improved capacity bounds
  for the binary energy harvesting channel,'' in {\em 2014 IEEE International
  Symposium on Information Theory}, pp.~976--980, IEEE, 2014.

\bibitem{DonEtal15}
Y.~Dong, F.~Farnia, and A.~{\"O}zg{\"u}r, ``Near optimal energy control and
  approximate capacity of energy harvesting communication,'' {\em IEEE Journal
  on Selected Areas in Communications}, vol.~33, no.~3, pp.~540--557, 2015.

\bibitem{OzeEtal15}
O.~Ozel, K.~Tutuncuoglu, S.~Ulukus, and A.~Yener, ``Fundamental limits of
  energy harvesting communications,'' {\em IEEE Communications Magazine},
  vol.~53, no.~4, pp.~126--132, 2015.

\bibitem{ShaEtal16}
D.~Shaviv, P.~Nguyen, and A.~{\"O}zg{\"u}r, ``Capacity of the energy-harvesting
  channel with a finite battery,'' {\em IEEE Transactions on Information
  Theory}, vol.~62, no.~11, pp.~6436--6458, 2016.

\bibitem{OzeUlu12b}
O.~Ozel and S.~Ulukus, ``On the capacity region of the {G}aussian {MAC} with
  batteryless energy harvesting transmitters,'' in {\em 2012 IEEE Global
  Communications Conference (GLOBECOM)}, pp.~2385--2390, IEEE, 2012.

\bibitem{InaEtal18}
H.~Inan, D.~Shaviv, and A.~{\"O}zg{\"u}r, ``Capacity of the energy harvesting
  {G}aussian {MAC},'' {\em IEEE Transactions on Information Theory}, vol.~64,
  no.~4, pp.~2347--2360, 2018.

\bibitem{OzeEtal11}
O.~Ozel, J.~Yang, and S.~Ulukus, ``Optimal scheduling over fading broadcast
  channels with an energy harvesting transmitter,'' in {\em 2011 4th IEEE
  International Workshop on Computational Advances in Multi-Sensor Adaptive
  Processing (CAMSAP)}, pp.~193--196, IEEE, 2011.

\bibitem{OzeEtal12}
O.~Ozel, J.~Yang, and S.~Ulukus, ``Optimal broadcast scheduling for an energy
  harvesting rechargeable transmitter with a finite capacity battery,'' {\em
  IEEE Transactions on Wireless Communications}, vol.~11, no.~6,
  pp.~2193--2203, 2012.

\bibitem{TutYen12}
K.~Tutuncuoglu and A.~Yener, ``Sum-rate optimal power policies for energy
  harvesting transmitters in an interference channel,'' {\em Journal of
  Communications and Networks}, vol.~14, no.~2, pp.~151--161, 2012.

\bibitem{YanUlu12}
J.~Yang and S.~Ulukus, ``Optimal packet scheduling in a multiple access channel
  with energy harvesting transmitters,'' {\em Journal of Communications and
  Networks}, vol.~14, no.~2, pp.~140--150, 2012.

\bibitem{ShaOzg16}
D.~Shaviv and A.~{\"O}zg{\"u}r, ``Universally near optimal online power control
  for energy harvesting nodes,'' {\em IEEE Journal on Selected Areas in
  Communications}, vol.~34, no.~12, pp.~3620--3631, 2016.

\bibitem{BakUlu16}
A.~Baknina and S.~Ulukus, ``Optimal and near-optimal online strategies for
  energy harvesting broadcast channels,'' {\em IEEE Journal on Selected Areas
  in Communications}, vol.~34, no.~12, pp.~3696--3708, 2016.

\bibitem{BakUlu18}
A.~Baknina and S.~Ulukus, ``Energy harvesting multiple access channels: Optimal
  and near-optimal online policies,'' {\em IEEE Transactions on
  Communications}, vol.~66, no.~7, pp.~2904--2917, 2018.

\bibitem{Sha48}
C.~E. Shannon, ``A mathematical theory of communication, {I} and {II},'' {\em
  Bell Syst. Tech. J}, vol.~27, pp.~379--423, 1948.

\bibitem{Smi71}
J.~G. Smith, ``The information capacity of amplitude-and variance-constrained
  scalar {G}aussian channels,'' {\em Information and Control}, vol.~18, no.~3,
  pp.~203--219, 1971.

\bibitem{ShaEtal15}
D.~Shaviv, A.~{\"O}zg{\"u}r, and H.~Permuter, ``Can feedback increase the
  capacity of the energy harvesting channel?,'' in {\em IEEE Information Theory
  Workshop (ITW), 2015}, pp.~1--5, IEEE, 2015.

\bibitem{Bla60}
D.~Blackwell, L.~Breiman, and A.~J. Thomasian, ``The capacities of certain
  channel classes under random coding,'' {\em The Annals of Mathematical
  Statistics}, vol.~31, no.~3, pp.~558--567, 1960.

\bibitem{CsiKor11}
I.~Csiszar and J.~K{\"o}rner, {\em Information theory: coding theorems for
  discrete memoryless systems}.
\newblock Cambridge University Press, 2011.

\bibitem{GolVar96}
A.~J. Goldsmith and P.~P. Varaiya, ``Capacity, mutual information, and coding
  for finite-state {M}arkov channels,'' {\em IEEE Transactions on Information
  Theory}, vol.~42, no.~3, pp.~868--886, 1996.

\bibitem{Wei10}
T.~Weissman, ``Capacity of channels with action-dependent states,'' {\em IEEE
  Transactions on Information Theory}, vol.~56, no.~11, pp.~5396--5411, 2010.

\bibitem{HirMas88}
W.~Hirt and J.~L. Massey, ``Capacity of the discrete-time {G}aussian channel
  with intersymbol interference,'' {\em IEEE Transactions on Information
  Theory}, vol.~34, no.~3, pp.~38--38, 1988.

\bibitem{Goldsmith2005Wireless-commun}
A.~Goldsmith, {\em Wireless communications}.
\newblock Cambridge University Press, 1~ed., 2005.

\bibitem{Tse2005Fundamentals-of}
D.~Tse and P.~Viswanath, {\em Fundamentals of Wireless Communication}.
\newblock Cambridge University Press, 1~ed., 2005.

\bibitem{Ozarow1994Information-the}
L.~H. Ozarow, S.~Shamai, and A.~D. Wyner, ``Information theoretic
  considerations for cellular mobile radio,'' {\em IEEE Transactions on
  Vehicular Technology}, vol.~43, no.~2, pp.~359--378, 1994.

\bibitem{Beaulieu2006A-closed-form-e}
N.~C. Beaulieu and J.~Hu, ``A closed-form expression for the outage probability
  of decode-and-forward relaying in dissimilar rayleigh fading channels,'' {\em
  IEEE Communications Letters}, vol.~10, pp.~813--815, December 2006.

\bibitem{Hadjiscostis2002Feedback-contro}
C.~N. Hadjiscostis and R.~Touri, ``Feedback control utilizing packet dropping
  links,'' in {\em Proceedings of the IEEE Conference on Decision and Control},
  pp.~1205--1210, 2002.

\bibitem{Sinopoli2004Kalman-filering}
B.~Sinopoli, L.~Schenato, M.~Franceschetti, K.~Poola, M.~I. Jordan, and S.~S.
  Sastry, ``Kalman filering with intermittent observations,'' {\em IEEE
  Transactions on Automatic Control}, vol.~49, no.~9, pp.~1453--1464, 2004.

\bibitem{Costa1993Stability-resul}
O.~L.~V. Costa and M.~D. Fragoso, ``Stability results for discrete-time linear
  systems with markovian jumping parameters,'' {\em Journal of Mathematical
  Analysis and Applications}, vol.~179, pp.~154--178, 1993.

\bibitem{Imer2006Optimal-control}
O.~C. Imer, S.~Yuksel, and T.~Basar, ``Optimal control of lti systems over
  unreliable communication links,'' {\em Automatica}, vol.~42, pp.~1429--1439,
  2006.

\bibitem{Imer2010Optimal-estimat}
O.~C. Imer and T.~Basar, ``Optimal estimation with limited measurements,'' {\em
  International Journal of Systems, Control and Communications}, vol.~2,
  pp.~5--29, 2010.

\bibitem{Schenato2007Foundations-of-}
L.~Schenato, B.~Sinopoli, M.~Franceschetti, K.~Poola, and S.~S. Sastry,
  ``Foundations of control and estimation over lossy networks,'' {\em
  Proceedings of the IEEE}, vol.~95, pp.~163--187, January 2007.

\bibitem{Xu2005Estimation-unde}
Y.~Xu and J.~P. Hespanha, ``Estimation under uncontrolled and controlled
  communications in networked control systems,'' in {\em Proceedings of the
  IEEE Conference on Decision and Control}, pp.~842--847, December 2005.

\bibitem{Rugh1996Linear-system-t}
W.~J. Rugh, {\em Linear system theory}.
\newblock Prentice Hall, 2~ed., 1996.

\bibitem{Hespanha2018Linear-systems-}
J.~P. Hespanha, {\em Linear systems theory}.
\newblock Princeton University Press, 2~ed., 2018.

\bibitem{Gupta2007Optimal-LQG-con}
V.~Gupta, B.~Hassibi, and R.~M. Murray, ``Optimal lqg control accross
  packet-dropping links,'' {\em Systems and Control Letters}, vol.~56,
  pp.~439--446, 2007.

\bibitem{Gupta2009Optimal-output-}
V.~Gupta, N.~C. Martins, and J.~S. Baras, ``Optimal output feedback control
  using two remote sensors over erasure channels,'' {\em IEEE Transactions on
  Automatic Control}, vol.~54, pp.~1463--1476, July 2009.

\bibitem{Gupta2010On-stability-in}
V.~Gupta and N.~C. Martins, ``On stability in the presence of analog erasure
  channel between the controller and the actuator,'' {\em IEEE Transactions on
  Automatic Control}, vol.~55, no.~1, pp.~175--179, 2010.

\bibitem{Minero2009Data-rate-theor}
P.~Minero, M.~Franceschetti, S.~Dey, and G.~N. Nair, ``Data rate theorem for
  stabilization over time-varying feedback channels,'' {\em IEEE Transactions
  on Automatic Control}, vol.~54, no.~2, pp.~243--255, 2009.

\bibitem{Wong1997Systems-with-Fi}
W.~S. Wong and R.~Brockett, ``Systems with finite communication bandwidth
  constraints---part i: State estimation problems,'' {\em IEEE Transactions on
  Automatic Control}, vol.~42, no.~9, 1997.

\bibitem{Martins2006Feedback-stabil}
N.~C. Martins, M.~A. Dahleh, and N.~Elia, ``Feedback stabilization of uncertain
  systems in the presence of a direct link,'' {\em IEEE Transactions on
  Automatic Control}, vol.~51, no.~3, pp.~438--447, 2006.

\bibitem{Sahai2006The-necessity-a}
A.~Sahai and S.~Mitter, ``The necessity and sufficiency of anytime capacity for
  stabilization of a linear system over a noisy communication link - part i:
  scalar systems,'' {\em IEEE Transactions on Information Theory}, vol.~52,
  pp.~3369--3395, August 2006.

\bibitem{Matveev2009Estimation-and-}
A.~S. Matveev and A.~Savkin, {\em Estimation and Control Over Communication
  Networks}.
\newblock Birkauser Boston, 2009.

\bibitem{Martins2011Remote-state-es}
G.~M. Lipsa and N.~C. Martins, ``Remote state estimation with communication
  costs for first-order lti systems,'' {\em IEEE Transactions on Automatic
  Control}, vol.~56, pp.~2013--2025, September 2011.

\bibitem{Hajek2008Paging-and-regi}
B.~Hajek, K.~Mitzel, and S.~Yang, ``Paging and registration in cellular
  networks: jointly optimal policies and an iterative algorithm,'' {\em IEEE
  Transactions on Information Theory}, vol.~54, pp.~608--622, February 2008.

\bibitem{Nayyar2013Optimal-strateg}
A.~Nayyar, T.~Basar, D.~Teneketzis, and V.~V. Veeravalli, ``Optimal strategies
  for communication and remote estimation with an energy harvesting sensor,''
  {\em IEEE Transactions on Automatic Control}, vol.~58, pp.~2246--2260,
  September 2013.

\bibitem{Park2018Individually-op}
S.~Park and N.~C. Martins, ``Individually optimal solutions to a remote state
  estimation problem with communication costs,'' in {\em Proceedings of the
  IEEE Conference on Decision and Control}, pp.~4014--4019, 2018.

\bibitem{Molin2013On-the-optimali}
A.~Molin and S.~Hirche, ``On the optimality of certainty equivalence for
  event-triggered control systems,'' {\em IEEE Transactions on Automatic
  Control}, vol.~58, no.~2, pp.~470--474, 2013.

\bibitem{Wang2011Event-triggerin}
X.~Wang and M.~Lemmon, ``Event-triggering in distributed networked control
  systems,'' {\em IEEE Transactions on Automatic Control}, vol.~56, pp.~586 --
  601, March 2011.

\bibitem{Puterman2005Markov-decision}
M.~L. Puterman, {\em Markov decision processes: discrete stochastic dynamic
  programming}.
\newblock Wiley, 2005.

\bibitem{Leong2018Transmission-sc}
A.~S. Leong, S.~Dey, and D.~E. Quevedo, ``Transmission scheduling for remote
  state estimation and control with an energy harvesting sensor,'' {\em
  Automatica}, vol.~91, pp.~54--60, 2018.

\bibitem{Li2017Power-control-o}
Y.~Li, F.~Zhang, D.~E. Quevedo, V.~Lau, S.~Dey, and L.~Shi, ``Power control of
  an energy harvesting sensor for remote state estimation,'' {\em IEEE
  Transactions on Automatic Control}, vol.~62, pp.~277--290, January 2017.

\bibitem{Nourian2014Optimal-Energy-}
M.~Nourian, A.~S. Leong, and S.~Dey, ``Optimal energy allocation for kalman
  filtering over packet dropping links with imperfect acknowledgements and
  energy harvesting constraints,'' {\em IEEE Transactions on Automatic
  Control}, vol.~59, pp.~2128--2143, August 2014.

\bibitem{Trimpe2014Event-based-sta}
S.~Trimpe and R.~D'Andrea, ``Event-based state estimation with variance-based
  triggering,'' {\em IEEE Transactions on Automatic Control}, vol.~59,
  pp.~3266--3281, December 2014.

\bibitem{Kumar2015Stochastic-Syst}
P.~R. Kumar and P.~Varayia, {\em Stochastic Systems: Estimation, Identification
  and Adaptive Control}.
\newblock SIAM, 2015.

\bibitem{Knorn2017Optimal-energy-}
S.~Knorn and S.~Dey, ``Optimal energy allocation for linear control with packet
  loss under energy harvesting constraints,'' {\em Automatica}, vol.~77,
  pp.~259--267, 2017.

\bibitem{Ren2016Infinite-Horizo}
X.~Ren, J.~Wu, K.~H. Johansson, G.~Shi, and L.~Shi, ``Infinite horizon optimal
  transmission power control for remote state estimation over fading
  channels,'' {\em IEEE Transactions on Automatic Control}, vol.~63,
  pp.~85--100, January 2018.

\bibitem{Chakravorty2019Remote-estimati}
J.~Chakravorty and A.~Mahajan, ``Remote estimation over a packet-drop channel
  with markovian state,'' {\em IEEE Transactions on Automatic Control (in
  press)}, 2019.

\bibitem{Lipsa2009Optimal-State-E}
G.~M. Lipsa and N.~C. Martins, ``Optimal state estimation in the presence of
  communication costs and packet drops,'' in {\em Proceedings of Allerton
  Conference on Communication, Control and Computing}, pp.~160--169, 2009.

\bibitem{Lin2018Remote-state-es}
M.~Lin, R.~J. La, and N.~C. Martins, ``Remote state estimation across an
  action-dependent packet-drop link,'' in {\em Proceedings of the IEEE
  Conference on Decision and Control}, 2018.

\bibitem{ConwayMaxwell62}
R.~W. Conway and W.~L. Maxwell, ``A queueing model with state dependent service
  rates,'' {\em Journal of Industrial Engineering}, vol.~12, pp.~132--136,
  1962.

\bibitem{Jackson63}
J.~R. Jackson, ``Jobshop-like queueing systems,'' {\em Management Science},
  vol.~10, no.~1, pp.~131--142, 1963.

\bibitem{Yadin63}
M.~Yadin and P.~Naor, ``Queueing systems with a removable service station,''
  {\em Operational Research Society}, vol.~14, pp.~393--405, December 1963.

\bibitem{Gupta67}
S.~Gupta, ``On bulk queues with state dependent parameters,'' {\em Journal of
  the Operations Research Society of Japan}, vol.~9, pp.~69--82, April 1967.

\bibitem{Harris67}
C.~M. Harris, ``Queues with state-dependent stochastic service rates,'' {\em
  Operations Research}, vol.~15, pp.~117--130, February 1967.

\bibitem{Dshalalow}
J.~H. Dshalalow, ``Queueing systems with state dependent parameters,'' in {\em
  Frontiers in Queueing: Models and Applications in Science and Engineering,
  Probability and Stochastics Series} (J.~H. Dshalalow, ed.), ch.~4,
  pp.~132--136, CRC, 1997.

\bibitem{YerkesDodson}
R.~M. Yerkes and J.~D. Dodson, ``The relation of strength of stimulus to
  rapidity of habit-formation,'' {\em Journal of Comparative Neurology and
  Psychology}, vol.~18, pp.~459--482, November 1908.

\bibitem{Edie54}
L.~C. Edie, ``Traffic delays at toll booths,'' {\em Journal of the Operations
  Research Society of America}, vol.~2, pp.~107--138, May 1954.

\bibitem{Borghini14}
G.~Borghini, L.~Astolfi, G.~Vecchiato, D.~Mattia, and F.~Babiloni, ``Measuring
  neurophysiological signals in aircraft pilots and car drivers for the
  assessment of mental workload, fatigue and drowsiness,'' {\em Neuroscience \&
  Biobehavioral Reviews}, vol.~44, pp.~58--75, July 2014.

\bibitem{Shunko17}
M.~Shunko, J.~Niederhoff, and Y.~Rosokha, ``Humans are not machines: the
  behavioral impact of queueing design on service time,'' {\em Management
  Science}, vol.~64, pp.~57--80, December 2017.

\bibitem{Sheridan97}
T.~S. Sheridan, ``Supervisory control,'' in {\em Handbook of Human Factors and
  Ergonomics, second edition} (G.~Salvendy, ed.), pp.~1295--1327, John Wiley \&
  Sons, 1997.

\bibitem{Asaro07}
P.~V. Asaro, L.~M. Lewis, and S.~B. Boxerman, ``The impact of input and output
  factors on emergency department throughput,'' {\em Adademic Emergency
  Medicine}, vol.~14, pp.~235--242, April 2007.

\bibitem{KcTerwiesch09}
D.~S. Kc and C.~Terwiesch, ``Impact of workload on service time and patient
  safety: an economic analysis of hospital operations,'' {\em Management
  Science}, vol.~55, pp.~1486--1498, September 2009.

\bibitem{Diamond07}
D.~M. Diamond, A.~M. Campbell, C.~R. Park, J.~Halonen, and P.~R. Zoladz, ``The
  temporal dynamics model of emotional memory processing: a synthesis on the
  neurobiological basis of stress-induced amnesia, flashbulb and traumatic
  memories, and the yerkes-dodson law,'' {\em Neural Plasticity}, 2017.

\bibitem{Lin-SQ}
M.~Lin, R.~J. La, and N.~C. Martins, ``Stability of a single queue subject to
  action-dependent server performance.'' preprint available at
  https://arxiv.org/abs/1903.00135, 2019.

\bibitem{Lin-Workload}
M.~Lin, N.~C. Martins, and R.~J. La, ``Queueing subject to action-dependent
  server performance: utilization rate reduction.'' preprint available at
  https://arxiv.org/abs/2002.08514, 2020.

\bibitem{Bi16}
S.~Bi, Y.~Zeng, and R.~Zhang, ``Wireless powered communication networks: An
  overview,'' {\em IEEE Wireless Communications}, vol.~23, pp.~10--18, April
  2016.

\bibitem{Niyato17}
D.~Niyato, D.~I. Kim, M.~Maso, and Z.~Han, ``Wireless powered communication
  networks: Research directions and technological approaches,'' {\em IEEE
  Wireless Communications}, vol.~24, pp.~88--97, December 2017.

\bibitem{Chu18}
Z.~Chu, F.~Zhou, Z.~Zhu, R.~Q. Hu, and P.~Xiao, ``Wireless powered sensor
  networks for internet of things: Maximum throughput and optimal power
  allocation,'' {\em IEEE Internet of Things Journal}, vol.~5, pp.~310--321,
  February 2018.

\bibitem{Lyu19}
B.~Lyu, T.~Qi, H.~Guo, and Z.~Yang, ``Throughput maximization in full-duplex
  dual-hop wireless powered communication networks,'' {\em IEEE Access},
  vol.~7, pp.~158584--158593, 2019.

\bibitem{Ju14}
H.~Ju and R.~Zhang, ``Throuhgput maximization in wireless powered communication
  networks,'' {\em IEEE Transactions on Wireless Communications}, vol.~13,
  pp.~418--428, January 2014.

\bibitem{Che15}
Y.~L. Che, L.~Duan, and R.~Zhang, ``Spatial throughput maximization of wireless
  powered communication networks,'' {\em IEEE Journal on Selected Ares in
  Communications}, vol.~33, pp.~1534--1548, August 2015.

\bibitem{Yang15}
G.~Yang, C.~K. Ho, R.~Zhang, and Y.~L. Guan, ``Throughput maximization for
  massive {MIMO} systems powered by wireless energy transfer,'' {\em IEEE
  Journal on Selected Ares in Communications}, vol.~33, pp.~1640--1650, August
  2015.

\bibitem{Shan16}
F.~Shan, J.~Luo, W.~Wu, and X.~Shen, ``Optimal wireless power transfer
  scheduling for delay minimization,'' in {\em Proceedings of the IEEE
  INFOCOM}, 2016.

\bibitem{Rezaei19}
R.~Rezaei, S.~Sum, X.~Kang, Y.~L. Guan, and M.~R. Pakravan, ``Secrecy
  throughput maximization for full-duplex wireless powered iot networks under
  fairness constraints,'' {\em IEEE Internet of Things Journal}, vol.~6,
  pp.~6964--6976, August 2019.

\bibitem{Kashef12}
M.~Kashef and A.~Ephremides, ``Optimal packet scheduling for energy harvesting
  sources on time varying wireless channels,'' {\em Journal of Communications
  and Networks}, vol.~14, pp.~121--129, April 2012.

\bibitem{Mao14}
S.~Mao, M.~H. Cheung, and V.~W. Wong, ``Joint energy allocation for sensing and
  transmission in rechargeable wireless sensor networks,'' {\em IEEE
  Transactions on Vehicular Technology}, vol.~63, pp.~2862--2875, July 2014.

\bibitem{Ahmed16}
I.~Ahmed, K.~T. Phan, and T.~Le-Ngoc, ``Optimal stochastic power control for
  energy harvesting systems with delay constraints,'' {\em IEEE Journal on
  Selected Areas In Communications}, vol.~34, pp.~3512--3527, December 2016.

\bibitem{AltmanMDP}
E.~Altman, {\em Constrained Markov Decision Processes}.
\newblock CRC Press, 1~ed., 1999.

\bibitem{Orhan2015Source-channel-}
O.~Orhan, D.~G\"{u}nd\"{u}z, and E.~Erkip, ``Source-channel coding under
  energy, delay and buffer constraints,'' {\em IEEE Transactions on Wireless
  Communcations}, vol.~14, no.~7, pp.~3836--3849, 2015.

\bibitem{Grover2011Towards-a-commu}
P.~Grover, K.~Woyach, and A.~Sahai, ``Towards a communication-theoretic
  understanding of system-level power consumption,'' {\em IEEE Journal on
  Selected Areas in Communications}, vol.~29, September 2011.

\bibitem{Wheeler2018Relative-Naviga}
D.~O. Wheeler, D.~P. Koch, J.~S. Jackson, T.~W. McLain, and R.~W. Beard,
  ``Relative navigation: a keyframe-based approach for observable gps-degraded
  navigation,'' {\em IEEE Control Systems Magazine}, vol.~38, pp.~30--48, July
  2018.

\bibitem{Kreucher2005Sensor-manageme}
C.~Kreucher, K.~Kastella, and A.~O. {Hero III}, ``Sensor management using an
  active sensing approach,'' {\em Signal Processing}, vol.~85, pp.~607--624,
  2005.

\bibitem{Walsh2001Scheduling-of-N}
G.~C. Walsh and H.~Ye, ``Scheduling of neyworked control systems,'' {\em IEEE
  Control Systems Magazine}, vol.~21, no.~1, pp.~57--65, 2001.

\bibitem{Sun2018Sampling-of-the}
Y.~Sun, Y.~Polyanskiy, and E.~Uysal{-}Biyikoglu, ``Sampling of the wiener
  process for remote estimation over a channel with random delay,'' {\em
  ArXiv}, 2018.

\bibitem{Sun2018Remote-Estimati}
Y.~Sun, Y.~Polyanskiy, and E.~Uysal{-}Biyikoglu, ``Remote estimation of the
  wiener process over a channel with random delay,'' {\em ArXiv}, 2018.

\bibitem{Sun2017Update-or-wait:}
Y.~Sun, E.~Uysal{-}Biyikoglu, R.~D. Yates, C.~E. Koksal, and N.~B. Shroff,
  ``Update or wait: how to keep your data fresh,'' {\em IEEE Transactions on
  Information Theory}, vol.~63, pp.~7492--7508, November 2017.

\bibitem{Arafa2018Age-minimal-tra}
A.~Arafa, J.~Yang, S.~Ulukus, and H.~V. Poor, ``Age-minimal transmission for
  energy harvesting sensors with finite batteries: online policies,'' {\em
  ArXiv}, 2018.

\bibitem{Kam2018Towards-an-effe}
C.~Kam, S.~Kompela, G.~D. Nguyen, J.~E. Wieselthier, and A.~Ephremides,
  ``{Towards an effective age of information: remote estimation of a Markov
  source},'' in {\em Proceedings of the IEEE Conference on Computer
  Communication Workshops: AoI Workshop}, pp.~367--372, 2018.

\bibitem{Ornee2019Sampling-for-re}
T.~Z. Ornee and Y.~Sun, ``Sampling for remote estimation through queues: age of
  information and beyond,'' {\em ArXiv}, February 2019.

\bibitem{Chatterjee17}
A.~Chatterjee, D.~Seo, and L.~R. Varshney, ``Capacity of systems with
  queue-length dependent service quality,'' {\em IEEE Transactions on
  Information Theory}, vol.~63, pp.~3950--3963, June 2017.

\bibitem{Pajic2017Design-and-impl}
M.~Pajic, J.~Weimer, N.~Bezzo, O.~Sokolsky, G.~J. Pappas, and I.~Lee, ``Design
  and implementation of attack-resilient cyberphysical systems,'' {\em IEEE
  Control Systems Magazine}, pp.~66--81, April 2017.

\bibitem{Fawzi2014Secure-estimati}
H.~Fawzi, P.~Tabuada, and S.~Diggavi, ``Secure estimation and control for
  cyber-physical systems under adversarial attacks,'' {\em IEEE Transactions on
  Automatic Control}, vol.~59, no.~6, pp.~1454--1467, 2014.

\bibitem{Shoukry2016Event-triggered}
Y.~Shoukry and P.~Tabuada, ``Event-triggered state observers for sparse sensor
  noise/attacks,'' {\em IEEE Transactions on Automatic Control}, vol.~61,
  pp.~2079--2091, August 2016.

\bibitem{Mo2010False-data-inje}
Y.~Mo, E.~Garone, A.~Casavola, and B.~Sinopoli, ``False data injection attacks
  against state estimation in wireless sensor networks,'' in {\em Proceedings
  of the IEEE Conference on Decision and Control}, pp.~5967--5972, 2010.

\bibitem{Jovanov2018Secure-state-es}
I.~Jovanov and M.~Pajic, ``Secure state estimation with cumulative message
  authentication,'' in {\em Proceedings of the IEEE Conference on Decision and
  Control}, pp.~2074--2079, 2018.

\bibitem{Sundaram2011Distributed-fun}
S.~Sundaram and C.~N. Hadjiscostis, ``Distributed function calculation via
  linear iterative strategies in the presence of malicious agents,'' {\em IEEE
  Transactions on Automatic Control}, vol.~56, pp.~1495--1508, July 2011.

\bibitem{Urbina2016Survey-and-new-}
D.~I. Urbina, J.~Giraldo, A.~A. Cardenas, J.~Valente, M.~Faisal, N.~O.
  Tippenhauer, J.~Ruths, R.~Candell, and H.~Sandberg, ``Survey and new
  directions for physics-based attack detection in control systems,'' Tech.
  Rep. NIST GCR 16-010, NIST, U.S. Department of Commerce, November 2016.

\bibitem{Cetinkaya2019An-overview-on-}
A.~Cetinkaya, H.~Ishii, and T.~Hayakawa, ``An overview on denial-of-service
  attacks in control systems: attack models and security analyses,'' {\em
  Entropy}, vol.~21, pp.~1--29, February 2019.

\bibitem{EasleyKleinberg}
D.~Easley and J.~Kleinberg, {\em Networks, Crowds and Markets}.
\newblock Cambridge University Press, 1~ed., 2010.

\bibitem{Dogar14}
F.~R. Dogar, T.~Karagiannis, H.~Ballani, and A.~Rowstron, ``Decentralized
  task-aware scheduling for data center networks,'' in {\em Proceedings of ACM
  SIGCOMM}, pp.~431--442, August 2014.

\bibitem{Mahajan2008Multi-armed-ban}
A.~Mahajan and D.~Teneketzis, {\em Multi-armed bandit problems},
  ch.~Multi-armed bandit problems, pp.~121--151.
\newblock Springer, 2008.

\bibitem{Whittle88}
P.~Whittle, ``Restless bandits: activity allocation in a changing world,'' {\em
  Journal of Applied Probability}, vol.~25, pp.~287--298, 1988.

\bibitem{RusZou16}
D.~Russo and J.~Zou, ``Controlling bias in adaptive data analysis using
  information theory,'' in {\em Proceedings of the 19th International
  Conference on Artificial Intelligence and Statistics} (A.~Gretton and C.~C.
  Robert, eds.), vol.~51 of {\em Proceedings of Machine Learning Research},
  pp.~1232--1240, PMLR, 09--11 May 2016.

\bibitem{XuRag17}
A.~Xu and M.~Raginsky, ``Information-theoretic analysis of generalization
  capability of learning algorithms,'' in {\em Advances in Neural Information
  Processing Systems} (I.~Guyon, U.~V. Luxburg, S.~Bengio, H.~Wallach,
  R.~Fergus, S.~Vishwanathan, and R.~Garnett, eds.), pp.~2521--2530, Curran
  Associates, Inc., 2017.

\bibitem{PenEtal18}
A.~Pensia, V.~Jog, and P.~Loh, ``Generalization error bounds for noisy,
  iterative algorithms,'' in {\em 2018 IEEE International Symposium on
  Information Theory (ISIT)}, pp.~546--550, June 2018.

\bibitem{NegEtal19}
J.~Negrea, M.~Haghifam, G.~K. Dziugaite, A.~Khisti, and D.~M. Roy,
  ``Information-theoretic generalization bounds for {SGLD} via data-dependent
  estimates,'' in {\em Advances in Neural Information Processing Systems},
  pp.~11013--11023, 2019.

\bibitem{SteEtal20}
T.~Steinke and L.~Zakynthinou, ``Reasoning about generalization via conditional
  mutual information,'' {\em arXiv preprint arXiv:2001.09122}, 2020.

\bibitem{GeEtal15}
R.~Ge, F.~Huang, C.~Jin, and Y.~Yuan, ``Escaping from saddle points---{O}nline
  stochastic gradient for tensor decomposition,'' in {\em Conference on
  Learning Theory}, pp.~797--842, 2015.

\bibitem{JinEtal17}
C.~Jin, R.~Ge, P.~Netrapalli, S.~M. Kakade, and M.~I. Jordan, ``How to escape
  saddle points efficiently,'' {\em arXiv preprint:1703.00887}, 2017.

\bibitem{AbbSan20}
E.~Abbe and C.~Sandon, ``Poly-time universality and limitations of deep
  learning,'' {\em arXiv preprint arXiv:2001.02992}, 2020.

\bibitem{BigEtal12}
B.~Biggio, B.~Nelson, and P.~Laskov, ``Poisoning attacks against support vector
  machines,'' in {\em Proceedings of the 29th International Conference on
  Machine Learning}, pp.~1467--1474, Omnipress, 2012.

\bibitem{MeiZhu15}
S.~Mei and X.~Zhu, ``Using machine teaching to identify optimal training-set
  attacks on machine learners.,'' in {\em AAAI}, pp.~2871--2877, 2015.

\bibitem{KohLia17}
P.~Koh and P.~Liang, ``Understanding black-box predictions via influence
  functions,'' {\em arXiv preprint arXiv:1703.04730}, 2017.

\bibitem{BlaEtal17}
P.~Blanchard, R.~Guerraoui, J.~Stainer, {\em et~al.}, ``Machine learning with
  adversaries: {B}yzantine tolerant gradient descent,'' in {\em Advances in
  Neural Information Processing Systems}, pp.~119--129, 2017.

\bibitem{CheEtal19}
Y.~Chen, L.~Su, and J.~Xu, ``Distributed statistical machine learning in
  adversarial settings: {B}yzantine gradient descent,'' {\em ACM SIGMETRICS
  Performance Evaluation Review}, vol.~46, no.~1, pp.~96--96, 2019.

\bibitem{Devillers2012A-General-Frame}
B.~Devillers and D.~G\"{u}nd\"{u}z, ``A general framework for the optimization
  of energy harvesting communication systems with battery imperfections,'' {\em
  Journal of Communications and Networks}, vol.~14, no.~2, pp.~130--139, 2012.

\bibitem{Tutuncouglu2015Optimum-Policie}
K.~Tutuncouglu, A.~Yener, and S.~Ulukus, ``Optimum policies for an energy
  harvesting transmitter under energy storage losses,'' {\em IEEE Journal on
  Selected Areas in Communications}, vol.~33, no.~3, pp.~467--481, 2015.

\end{thebibliography}

\end{document}